\documentclass{amsart}
\usepackage[utf8]{inputenc}

\usepackage{amssymb}
\usepackage{amsfonts}
\usepackage{amsmath}
\usepackage{amstext}
\usepackage{amsthm}
\usepackage{latexsym}

\usepackage[utf8]{inputenc}

\usepackage{mathrsfs}  
\usepackage[mathscr]{eucal}

\usepackage{amscd}
\usepackage{eucal}
\usepackage{amsxtra}
\usepackage{amsbsy}
\usepackage{graphicx}
\usepackage{latexsym}
\usepackage{pb-diagram}
\usepackage{amsopn}
\usepackage{delarray}
\usepackage{psfrag}
\usepackage{verbatim}
\usepackage{bussproofs}
\usepackage[T1]{fontenc}
\usepackage{color,cancel}

\usepackage{appendix}
\usepackage{color}
\usepackage[dvipsnames]{xcolor}
\usepackage{epigraph}
\usepackage{lipsum}
\usepackage{hyperref} %LaTex Hyperlinks
\hypersetup{pdfpagemode={UseOutlines},
bookmarksopen=true,bookmarksopenlevel=0,hypertexnames=false,colorlinks=true,citecolor=black,linkcolor=black,urlcolor=black,pdfstartview={FitV},unicode,breaklinks=true}
\usepackage{graphicx}
\usepackage{enumitem}

\newtheorem{thm}{Theorem}[section]
\newtheorem{propo}[thm]{Proposition}
\newtheorem{coro}[thm]{Corollary}
\newtheorem{lem}[thm]{Lemma}
\newtheorem{que}[thm]{Question}
\newtheorem*{que*}{Question}
\newtheorem{con}[thm]{Conjecture}

\theoremstyle{definition}
\newtheorem{dfn}[thm]{Definition}

\theoremstyle{remark}

\newtheorem{fact}[thm]{Fact}
\newtheorem{nota}[thm]{Notation}

\newtheorem{example}[thm]{Example}

\newcommand{\conc}{^\smallfrown}	% Concatenation.
\newcommand{\rest}{\upharpoonright} % Restriction.
\newcommand{\QED}{\hfill\ensuremath{\square}}

\def\parr{{\mathfrak{par}}}
\def\homm{{\mathfrak{hom}}}

\def\HH{{\mathbb{H}}}

\def\cA{{\mathcal{A}}} \def\cB{{\mathcal{B}}} \def\cC{{\mathcal{C}}} \def\cD{{\mathcal{D}}}  \def\cF{{\mathcal{F}}} \def\cG{{\mathcal{G}}} \def\cH{{\mathcal{H}}} \def\cI{{\mathcal{I}}} \def\cJ{{\mathcal{J}}}    \def\cN{{\mathcal{N}}}  \def\cP{{\mathcal{P}}}  \def\cR{{\mathcal{R}}} \def\cS{{\mathcal{S}}}  \def\cV{{\mathcal{V}}} \def\cU{{\mathcal{U}}}  \def\cX{{\mathcal{X}}} \def\cY{{\mathcal{Y}}} \def\cZ{{\mathcal{Z}}}

                \def\bQ{{\mathbb{Q}}} \def\bR{{\mathbb{R}}}        

\def\a{{\mathfrak{a}}}
\def\b{{\mathfrak{b}}}
\def\c{{\mathfrak{c}}}
\def\d{{\mathfrak{d}}}

\def\p{{\mathfrak{p}}}

\def\r{{\mathfrak{r}}}
\def\s{{\mathfrak{s}}}

\def\u{{\mathfrak{u}}}

\def\fin{{\mathsf{FIN}}}%we have this to write the fin ideal

\newcommand{\w}{\omega}
\newcommand{\Katetov}{Kat\v{e}tov }
\newcommand{\seqcomp}{sequentially compact }
\newcommand{\seqcompp}{sequentially compact}%the same as the later but without the space at the end, this can be used when a phrase finishes with "sequentially copmpact"

\title[Infinite dimensional sequential compactness]{Infinite dimensional sequential compactness: Sequential compactness based on barriers}

\author[Corral]{C. Corral}
\address{York University, 4700 Keele St, Toronto, ON, Canada, M3J 1P3.}
\email{cicorral@yorku.ca}
\thanks{The first author acknowledges support from York University and the Fields Institute.}

\author[Guzm\'an]{O. Guzm\'an}
\address{Centro de Ciencias  Matem\'aticas, Universidad Nacional Aut\'onoma de M\'exico, Campus Morelia, 58089, Morelia, Michoac\'an, M\'exico.}
\email{oguzman@matmor.unam.mx}
\thanks{The research of the second author was supported by PAPIIT grant IA102222.}

\author[L\'opez]{C. L\'opez-Callejas}
\address{Centro de Ciencias  Matem\'aticas, Universidad Nacional Aut\'onoma de M\'exico, Campus Morelia, 58089, Morelia, Michoac\'an, M\'exico.}
\email{carloscallejas@matmor.unam.mx}
\thanks{The research of the third author was supported by PAPIIT grant IN101323 and CONACyT grant A1-S-16164.}

\author[Memarpanahi]{P. Memarpanahi}
\address{University of Toronto Scarborough 1265 Military Trail, Toronto, Ontario M1C 1A4, Canada}
\email{pourya.memarpanahi@utoronto.ca}

\author[Szeptycki]{P. Szeptycki}
\address{Department of Mathematics and Statistics, York University, Toronto, Ontario M3J 1P3, Canada}
\email{szeptyck@yorku.ca}
\thanks{The fifth author acknowledges support from NSERC}

\author[Todor\v{c}evi\'c]{S. Todor\v{c}evi\'c}
\address{Department of Mathematics, University of Toronto, Canada}
\email{stevo@math.toronto.edu}
\address{Institut de Math\'ematiques de Jussieu, CNRS, Paris, France}
\email{stevo.todorcevic@imj-prg.fr}

\address{Matemati\v{c}ki Institut, SANU, Belgrade, Serbia}
\email{stevo.todorcevic@sanu.ac.rs}
\thanks{The research of the sixth author was partially supported by grants from
NSERC(455916), CNRS(UMR7586), SFRS(7750027-SMART) and EXPRO 20-31529X
grant (Czech Science Foundation)}

\address{Mathematical Institute of the Czech Academy of Sciences}

\date{}



\begin{document}

\maketitle

\begin{abstract}
We introduce a generalization of sequential compactness using barriers on $\w$ extending naturally the notion introduced in \cite{kubistopologicalramsey}. We improve results from \cite{corralnRamsey} by building spaces that are $\cB$-sequentially compact but no $\cC$-sequentially compact when the barriers $\cB$ and $\cC$ satisfy certain rank assumption which turns out to be equivalent to a Kat\v{e}tov-order assumption. Such examples are constructed under the assumption  $\b =\c$. We also exhibit some classes of spaces that are $\cB$-sequentially compact for every barrier $\cB$, including some classical classes of compact spaces from functional analysis, and as a byproduct we obtain some results on angelic spaces. Finally we introduce and compute some cardinal invariants naturally associated to barriers.
\end{abstract}

%\tableofcontents

\section{Introduction}

A 2-dimensional version of sequential compactness 
%that we generalize here to infinite dimensions 
was first considered by M. Boja\'nczyk, E. Kopczy\'nski and S. Toru\'nczyk  in \cite{ramseymetric} where they showed that if $f:[\w]^2\to K$ and $K$ is a compact metric space, then there is an infinite set $B\in[\w]^\w$ and $x\in K$ such that for every open $U\ni x$, there exists $m\in\w$ such that $f''[B\setminus m]^2\subseteq U$. In this case, $x$ is said to be the limit of $f\rest [B]^2$.
Using this result, they show that every compact metric semigroup has an idempotent that can be defined as the limit of naturally defined $f:[\w]^2\to K$. It is natural
to call this property $2$-dimensional sequential compactness and look for higher-dimensional versions. In fact, the notion of $n$-sequential compactness was introduced and studied in \cite{kubistopologicalramsey} (called $n$-Ramsey in their paper) as a natural generalization from the case $n=2$ to all other positive integers $n$. A space $X$ is said to be \emph{$n$-sequentially compact} if for every function $f:[\w]^n\to X$ there is an infinite set $M\in[\w]^\w$ and $x\in X$ such that for every open $U\ni x$ there exists $k\in\w$ such that $f''[M\setminus k]^n\subseteq U$.
It was shown in \cite{kubistopologicalramsey} that every $n$-sequentially compact space is also $m$-sequentially compact as long as $m\leq n$, and since $1$-sequential compactness coincides with classical notion of sequential compactness, every $n$-sequentially compact space is also sequentially compact. On the other hand, they gave examples of $n$-sequentially compact spaces that are not $(n+1)$-sequentially compact (assuming \textsf{CH} for $n>1$) and proved that sequentially compact spaces of character less that $\b$ are $n$-sequentially compact for all $n\in\w$ among other results.
Many of these results were further improved in \cite{corralnRamsey}, where the study of $n$-sequentially compact spaces was carried on.

\vspace{2mm}

It should be mentioned that a related notion of {\it Ramsey convergence} was introduced and studied by H. Knaust in \cite{knaustarray}. Given an array $\{x_{i,j}:0\leq i<j<\w\}$ in a space $X$, it is said that it converges \emph{Ramsey-uniformly} to a point $x$ in $X$ if there is an infinite set $M\in[\w]^\w$ such that for every open neighborhood $U$ of $x$, there is $k\in\w$ such that $x_{i,j}\in U$ for every $i,j\in M$ with $i<j$. 
Knaust also defines a space $X$ to have the {\it Ramsey property} if given an array $\{x_{n,m}:n<m<\w\}$ and a point $x$ in the space $X$ such that $\lim_{n\to\infty}\lim_{m\to\infty}x_{n,m}=x$, the array $\{x_{n,m}:n<m<\w\}$ converges Ramsey uniformly to $x$.
He then showed that every Rosenthal compact has the Ramsey property. 
In  the  subsequent paper \cite{knaustangelic}, Knaust showed that some classes of angelic spaces have the Ramsey property, including function spaces $C_p(X)$ over quasi-Suslin spaces $X$.

%The notion of an $n$-sequentially compact spaces was introduced in \cite{kubistopologicalramsey} with different terminology, and studied further in \cite{corralnRamsey} using the current terminology. We generalize this notion to general barriers in $[\omega]^{<\omega}$ (see the definitions in Section \ref{barrierssection}) and consider the stronger properties of a space being $\cB$-sequentially compact or more generally $\alpha$-sequentially compact where $\alpha$ denotes the rank of a given barrier. 
The  purpose of this paper is to extend the notion of sequential compactness to infinite dimensions. This is done using {\it barriers} of Nash-Williams \cite{NashWilliams}, where the notion of ``dimension'' is captured by the rank of the barrier which can be any countable ordinal. Since the simplest examples of barriers are the families $[\omega]^n$, the notion of ``Barrier-sequential compactness'' is a natural one. 

%notion of  \cite{ramseymetric} and  \cite{kubistopologicalramsey} which in hindsight use the barriers $[\omega]^n$ for $1\leq n<\omega$ of finite ranks.

We start in Section \ref{barrierssection} by studying the general theory of barriers and stating many results that are the main tools used in our proofs in later sections.

In Section \ref{sectionpositiveresults} we study the classes of $\alpha$-sequentially compact spaces and prove some basic results. For this definition to make sense, we will show that this property depends mostly on the rank of a given barrier and not on the particular recursive structure of it. 

In Section \ref{examples} we show that sequentially compact spaces of small character and compact bisequential spaces are $\alpha$-sequentially compact for every $\alpha<\omega_1$. As a Corollary, we obtain that many classical classes of compact spaces are $\alpha$-sequentially compact for every $\alpha$. We close the section with some comments on angelic spaces and we point out that a space constructed in \cite{kubistopologicalramsey} and improved in \cite{corralnRamsey} is an angelic space that fails to satisfy the Ramsey property, answering a question of Knaust \cite{knaustarray}. 

In Section \ref{counterexamplessection} we present the constructions of spaces that are $\beta$-sequentially compact but not $\alpha$-sequentially compact for $\beta<\alpha$, under the assumption of $\mathfrak{b}=\mathfrak{c}$. We also analyze further the properties mentioned in Section \ref{sectionpositiveresults} and their relation to $\alpha$-sequential compactness.

Finally, in section \ref{Cardinalinvariants} we define and discuss a number of cardinal invariants associated to these classes of spaces, and give some topological applications.

Our terminology and notation is mostly standard. In particular $\cN(x)$ stands for the collection of open neighborhoods of $x$ and for a finite set $\{b_0,\ldots,b_n\}\in[\w]^{<\w}$, we will always assume that it is written in increasing order.

We will be using a number of classical cardinal invariants of the continuum including $\p\a\r$, $\p$ $\b$, $\s$, $\d$, $\r$, $\r_\sigma$. We refer the reader to \cite{blasscardinal} for the background on these cardinals.

%Our counterexamples are almost exclusively compactifications of $\Psi$-spaces constructed over almost disjoint families of infinite subsets of $\omega$. For more on ad families and $\Psi$ spaces see \cite{michelferpsispaces}. 

Given an almost disjoint family $\cA$ on a countable set $N$, we define its \emph{Franklin space} as then the one point compactification of $\Psi(\cA)$, i.e., $\cF(\cA)=\Psi(\cA)\cup\{\infty\}$. Our counterexamples are almost exclusively Franklin spaces. For more on almost disjoint families and $\Psi$-spaces see \cite{michelferpsispaces}. 

Any other set-theoretic notion and terminology is standard and can be found in \cite{kunenbook} and we refer the reader to \cite{engelking} for topological terminology.

%At some point, may be in the introduction or in a further section, we have to talk about ultrafilters on barriers.\cesar{I will be using the tree version for these ultrafilters, but maybe we should state all the equivalent forms and use the most appropriate for each purpose}.

%$\cF(\cA)=\Psi(\cA)^*$ is the Franklin space of the ad family $\cA$.

%$\cA$ satisfies $P$ if and only if its Franklin space.

%We identify finite sets $s$ with its increasing enumeration $s^\uparrow$ and families of finite sets $\cC$ with the set $\cC^\uparrow=\{s^\uparrow:s\in\cC\}$. Imagenes de funciones, mencionar ambas notaciones.

%Trees.

%Notation $S^0=S$ and $S^1=\w\setminus S$.

%Antichain definition.

\section{Barriers and fronts}\label{barrierssection}

We now introduce the basic concepts of barriers and fronts on $[\omega]^{<\omega}$.  This section is mainly based on \cite{Todorcevicramseymethods} and \cite{todorcevicramseyspaces}, although some results are new or adapted and will be applied in subsequent sections where we will extend the theory of $n$-sequentially compact spaces by defining the notion of barrier sequential compactness and ultimately $\alpha$-sequential compactness for all $\alpha<\omega_1$.

For $\cF\subseteq[\omega]^{<\omega}$ and $M\in[\omega]^\omega$ let 
$$\cF|M=\{s\in\cF:s\subseteq M\}.$$
Recall that $s\sqsubseteq t$ means that $t$ is an end extension of $s$ and $s\sqsubset t$ means $s\sqsubseteq t$ and $s\neq t$.

\begin{dfn}
    A family $\cF\subseteq[\omega]^{<\omega}$ is:
    \begin{itemize}
        \item \emph{Ramsey}: if for every partition $\cF=\cF_0\sqcup\cdots\sqcup\cF_n$ and for every $N\in[\omega]^\omega$, there is $M\in[N]^\omega$ such that all but at most one of the restrictions $\cF_i|M$ are empty.
        \item \emph{Nash-Williams:} if $s\sqsubseteq t$ implies $s=t$ for $s,t\in\cF$.
        \item \emph{Sperner:} if $s\subseteq t$ implies $s=t$ for $s,t\in \cF$.
    \end{itemize}
\end{dfn}

Given a Ramsey family $\cF$, applying the Ramsey property to the partition $\cF=\cF_0\cup(\cF\setminus\cF_0)$ where $\cF_0$ is the set of $\subseteq$-minimal elements in $\cF$ we get the following:

\begin{propo}\label{Ramseyisspermerpropo}
    If $\cF$ is Ramsey, there exists $M\in[\omega]^\omega$ such that $\cF|M$ is Sperner.
\end{propo}

\begin{thm}\cite{NashWilliams}\label{NashWilliamsareRamsey}
    Every Nash-Williams family is Ramsey.
\end{thm}

The previous two results show that (at least if one is willing to pass to an infinite subset) the three concepts of being Ramsey, Nash-Williams and Sperner are equivalent:
\[\textnormal{Ramsey}\Rightarrow\textnormal{(some restriction is) Sperner}\Rightarrow\textnormal{Nash-Williams}\Rightarrow\textnormal{Ramsey.}\]

\begin{dfn}
A Nash-Williams family $\cF$ such that for every $M\in[\omega]^\omega$ there is an initial segment $s\sqsubseteq M$ with $s\in\cF$ is called a \emph{front}.
If moreover $\cF$ is Sperner, we say that $\cF$ is a \emph{barrier} on $\omega$.
\end{dfn}

The same argument for Proposition \ref{Ramseyisspermerpropo} can be used to show the following facts:

\begin{fact} If $\cF$ is a front on a countable set $M$, there exists $N\in[M]^\omega$ such that $\cF|N$ is a barrier.
\end{fact}

\begin{fact}
For every barrier $\cB$ on $M$, if we partition $\cB=\cB_0\sqcup\cdots\sqcup\cB_n$, there are $N\in[M]^\omega$ and $i\leq n$ such that $\cB_i|N$ is a barrier.
\end{fact}

Given a barrier $\cB$, let 
$$T(\cB)=\{s\in[\omega]^{<\omega}:\exists t\in\cB\ (s\sqsubseteq t)\}$$ 
and $\rho_{T(\cB)}:T(\cB)\to\omega_1$ be given by
$$\rho_{T(\cB)}(s)=\sup\{\rho_{T(\cB)}(t)+1:t\in T(\cB)\land s\sqsubset t\},$$
where $\sup\emptyset=0$. We will omit the subindex $T(\cB)$ or replace it by $\cB$ when no confusion arises and we will often think of $T(\cB)$ as a subtree of $\w^{<\w}$.

\begin{dfn}
    For a barrier $\cB$ on a countable set $M$, its rank is defined by $\rho(\cB):=\rho_{T(\cB)}(\emptyset)$.
\end{dfn}

It is worth to point out that in the definition of $T(\cB)$, we can replace the condition of being an initial segment of an element $t\in\cB$, to being a subset of some $t'\in\cB$ due to the following fact.

\begin{lem}\label{contentionandinitialsegment}
    If $\cB$ is a barrier on $\omega$, then for every $s\in[\w]^{<\w}$, there exists $b\in \cB$ such that $s\subseteq b$ if and only if there exists $b\in \cB$ such that $s\sqsubseteq b$. In particular $T(\cB)=\{s\in[\w]^{<\w}:\exists b\in\cB\ (s\subseteq b)\}$.
\end{lem}

\begin{proof}
    Let $s\in[\w]^{<\w}$ and $b\in\cB$ such that $s\subseteq b$. Define $M=s\cup(\w\setminus\max(b)+1)$. We can find $b'\in \cB$ such that $b'\sqsubseteq M$ and hence $s\sqsubseteq b'$ since otherwise $b'\sqsubset s\subseteq b$ would contradict that $\cB$ is a $\subseteq$-antichain.
\end{proof}

Given $a=\{a_0,\ldots,a_n\}$ and $b=\{b_0,\ldots,b_k\}$ with $k<n$, denote by $a\ast b$ the \emph{end replacement} of $a$ with $b$, defined as follows:
\[a\ast b=\{a_0,\ldots,a_{n-k-1},b_0,\ldots,b_k\}.\]

\begin{lem}\label{homogeneityofbarriers}
    Let $B$ be a barrier and $a\in\cB$ enumerated as $a=\{a_0,\ldots,a_{n-1}\}$. Then for every $b\in[\omega\setminus(a_{n-1})]^{<n}$, there exists $s\in\cB$ such that $a\ast b\sqsubseteq s$. In particular, $b\notin\cB$. 
\end{lem}

\begin{proof}
Fix $a=\{a_0,\ldots,a_{n-1}\}$ and $b=\{b_0,\ldots,b_{k-1}\}$ as in the lemma with $k<n$. We will prove the Lemma by induction on $k$. 

If $k=1$ then $b=\{b_0\}$ for some $b_0\geq a_{n-1}$. Let $M\in[\omega]^\omega$ with $a\ast b\sqsubseteq M$. 
There is an $s\in\cB$ such that $s\sqsubseteq M$. Notice that $|s|\geq n$ since otherwise $s\sqsubset a$ would contradict that $\cB$ is Sperner. Thus $a\ast b\sqsubseteq s\in\cB$.

Now let $k>1$ and assume the result is true for ${k-1}$. Pick any $M\in[\omega]^\omega$such that $a\ast b\sqsubseteq M$ and let $s\in\cB$ be an initial segment of $M$. If $|s|<n$, then $s\sqsubseteq\{a_0,\ldots,a_{n-k-1},b_0,\ldots,b_{k-2}\}\subseteq\{a_0,\ldots,a_{n-k},b_0,\ldots,b_{k-2}\}=a\ast b'$ where $b'=\{b_0,\ldots,b_{k-2}\}$. Since $|b'|=k-1$, by the inductive hypothesis we get that $s\subsetneq a\ast b'\sqsubseteq s'$ for some $s'\in\cB$, contradicting that $\cB$ is Sperner. Therefore $|s|\geq n$ and hence $a\ast b\sqsubseteq s\in\cB$.
\end{proof}

\begin{coro}\label{uniquenessoffiniterankbarriers}
If $\cB$ is a barrier on $\w$ of rank $k\in\mathbb{N}$, there exists $m\in\omega$ such that $\cB|(\omega\setminus m)=[\omega\setminus m]^k$.
\end{coro}

\begin{proof}
    Let $\cB$ be a barrier of rank $k$. If there is $a\in\cB$ of size $l$, then $\rho(\cB)\geq l$. Thus $\cB\subseteq[\omega]^{\leq k}$. It is also easy to see that $\cB\cap[\omega]^k\neq \emptyset$ since otherwise $\rho(\cB)<k$.

    Take any $a\in\cB\cap[\omega]^k$ and let $m=\max(a)+1$. By Lemma \ref{homogeneityofbarriers}, if $b\in\cB|(\w\setminus m)$ has size less than $k$, we can find $s\in\cB$ such that $b'=a\ast b\sqsubseteq s$. As $|b|<k$ we have that $b\subsetneq b'\sqsubseteq s$, which contradicts that $\cB$ is Sperner. Henceforth $\cB|(\w\setminus m)=[\w\setminus m]^k$.
\end{proof}

\begin{nota}
Given a family $\cB\subseteq[\w]^{<\omega}$, $s\in[\w]^{<\w}$ and $n\in\omega$ set
\begin{itemize}
    \item $\cB(s)=\{t\in\cB:s\sqsubseteq t\}$,
    \item $\cB[s]=\{t\setminus s:t\in\cB(s)\}$,
    \item $s\conc\cB=\{s\cup t:t\in\cB\ \land\ \min(t)>\max(s)\}$,
    \item $\cB+n=\{s+n:s\in\cB\}$, where $s+n=\{m+n:m\in s\}$.
\end{itemize}
\end{nota}

We write $\cB(n)$, $\cB[n]$ and $n\conc\cB$ instead of $\cB(\{n\})$, $\cB[\{n\}]$ and $\{n\}\conc\cB$ respectively. Notice that some confusion may arise as $n=\{0,\ldots,n-1\}$ but this notation won't lead to any confusion as $\cB(n)$ will always be interpreted as $\cB(\{n\})$ and the same for $\cB[n]$ and $n\conc\cB$.
We see that if $\cB$ is a barrier on $M$, then $\cB[n]$ is a barrier on $M\setminus(n+1)$.
Conversely, if $\cB_n$ is a barrier on $M\setminus (n+1)$ for every $n\in M$, then $\bigcup_{n\in M}n\conc\cB_n$ is a front on $M$ and there is an infinite set on which its restriction is a barrier. We can now describe a canonical barrier of rank $\w$: The \emph{Schreier barrier} $\cS$ is the barrier of rank $\w$ such that $\{0\}\in\cS$ and $\cS[n]=[\w\setminus(n+1)]^n$ for every $n>0$. In other words, $s\in\cS$ if and only if $|s|=\min(s)+1$.

\begin{dfn}
    Let $\cB$ be a barrier on $M\in[\w]^\w$ with $\rho(\cB)=\alpha$. If $\alpha>1$, we say that $\cB$ is an \emph{uniform} barrier if each $\cB[n]$ is an uniform barrier (on $M\setminus(n+1)$) and
    \begin{itemize}
        \item $\rho_{T(\cB)}(\{n\})=\beta$ for every $n\in\w$ if $\alpha=\beta+1$ or
        \item $\{\rho_{T(\cB)}(\{n\}):n\in\omega\}$ is an increasing sequence with limit $\alpha$ if $\alpha$ is limit.
    \end{itemize}
    For $\alpha=1$, we declare the unique barrier $\cB=[\w]^1$ as an uniform barrier. 
\end{dfn}

The following result together with Corollary \ref{uniquenessoffiniterankbarriers} shows that all barriers of rank $\omega$ are somewhere uniform and preserve its rank in such restriction.

\begin{propo}\label{Allbarriersofrankomegaareuniform}
If $\cB$ is a barrier of rank $\omega$, there exists $M\in[\omega]^\omega$ such that $\cB|M$ is uniform and has rank $\omega$.
\end{propo}

\begin{proof}
    By definition, we can find $m_i$ such that $\rho(\{m_i\})\geq i$ for every $i\in\omega$. By Corollary \ref{uniquenessoffiniterankbarriers}, we can also find $k_i$ for every $i\in\omega$, such that $\cB[{m_i}]=[\omega\setminus k_i]^{n_i}$, where $n_i=\rho(\{m_i\})$. Notice that we can pick $\{m_i:i\in\omega\}$ increasing and such that $m_{i+1}>k_i$. It is now easy to see that $M=\{m_i:i\in\omega\}$ works.
\end{proof}

A partial analogue for Proposition \ref{Allbarriersofrankomegaareuniform} for any barrier independently of its rank is also true.

\begin{propo}\cite{Todorcevicramseymethods}
    For every barrier $\cB$ on $\w$ there exists an infinite set $M\in[\w]^\w$ such that $\cB|M$ is uniform.\QED
\end{propo}

By Corollary \ref{uniquenessoffiniterankbarriers}, the only uniform barrier of rank $k\in\omega$ is $[\omega]^k$. Thus Lemma \ref{homogeneityofbarriers} states that the sequence of ranks $\{\rho(\{n\}):n\in\omega\}$ is non decreasing and the sequence $\{m_i:i\in\omega\}$ can be taken to be the identity in the previous theorem. 
This shows that the only way to get a non-uniform barrier of rank $\omega$ is by embedding a non-uniform barrier $\cB[n]$ on top of $\{n\}$, hence essentially, every barrier of rank $\omega$ is uniform. This is as far as we can go since there is a non uniform barrier $\cB$ of rank $\w+1$ such that $\cB[n]$ is uniform for every $n\in\w$.

\begin{example}
    There is a non-uniform barrier of rank $\w+1$.

    \vspace{1mm}
    
    We define a barrier $\cB$ by describing $\cB(n)$ for every $n\in\w$. For every $k\in\w$, let $\cS_k=\{s\in[\w]^{<\w}:|s|=\min(s)+k\}$, (hence the Schreier barrier is $\cS_1$). 
    If $n$ is even define $\cB(n)=n\conc(\cS_n+n+1)$. For $n$ odd let $\cB(n)=n\conc[\omega\setminus(n+1)]^n$. 

    It is clear that every infinite subset of $\w$ has an initial segment in $\cB$ and that each $\cB(n)$ is a $\subseteq$-antichain. 
    It remains to show that if $s=\{s_0,\ldots,s_i\}\in\cB(n)$, $t=\{t_0,\ldots,t_j\}\in\cB(m)$ and $n<m$, then $t\nsubseteq s$ (the other contention is impossible as $n\in s\setminus t$).

    Let us first compute the size of an element $b\in\cB$. Let $b=\{b_0,\ldots,b_k\}$. If $b_0$ is odd, then $b\in\cB(b_0)$ and it has size $b_0+1$. Otherwise, $b\in b_0\conc(\cS_{b_0}+b_0+1)$. Let $b'=b\setminus\{b_0\}\in \cS_{b_0}+b_0+1$. It is clear that $|b'|=|b'-b_0-1|$ and since $(b'-b_0-1)\in\cS_{b_0}$ and $\min(b'-b_0-1)=b_1-b_0-1$, we conclude that $|b'|=|b'-b_0-1|=(b_1-b_0-1)+b_0=b_1-1$. Therefore $|b|=b_1$.

    We are now ready to show that $s$ and $t$ are $\subseteq$-incomparable.
    If both $n$ and $m$ are odd we have that $|t|=m+1>n+1=|s|$ and we are done.
    If both $n$ and $m$ are even, then $|s|=s_1$ and $|t|=t_1$ but $s_0=n\not\in t$ implies that if 
    $t\subseteq s$, hence $|s|=s_1\leq t_0<t_1=|t|$, which is a contradiction.\\
    Similarly if $n$ is odd and $m$ is even, assuming that $t\subseteq s$ yields that $|t|=t_1>t_0\geq s_1\geq s_0+1=|s|$, a contradiction.
    Finally, if $n$ is even and $m$ is odd, and if we assume that $t\subseteq s$ we reach a contradiction as $|t|=t_0+1>t_0\geq s_1=|s|$.
\QED
\end{example}

The persistence of the rank under restriction for a uniform barrier is the key property of the notion of uniformity. We can imitate this behaviour by considering a lower bound for the ranks of uniform restrictions of a given barrier $\cB$.

\begin{dfn}
Given a barrier $\cB$, let 
$$spec(\cB)=\{\alpha<\omega_1:\exists M\in[\w]^\w\ (\cB|M\textnormal{ is uniform with rank } \alpha)\}.$$

We define the \emph{uniform rank of} $\cB$ as $\rho_u(\cB)=\min(spec(\cB))$.
\end{dfn}

%\vspace{3mm}

In \cite{kubistopologicalramsey}, it is inductively proved  that sequentially compact spaces of character less that $\b$ are $n$-sequentially compact for every $n$. To perform the induction on $[\w]^{n+1}$, it is enough to notice that any element of $[\w]^{n+1}$ has a unique initial segment in $[\w]^n$. It is also true that any element in $[\w]^n$ is end-extended by at least one element in $[\w]^{n+1}$. This fact is used to analyze a splitting-like cardinal invariant in \cite{corralnRamsey} that helps with the construction of $n$-sequentially compact spaces that fail to be $(n+1)$-sequentially compact under some assumptions involving this cardinal. Other examples of results that use this fact, are Ramsey's theorem itself and the proof that $\p\a\r_n=\p\a\r_2=\max\{\b,\s\}$ in \cite{blasscardinal}. In order to obtain analogous results in our framework, we will need a generalization of this fact for barriers of any rank. The statement of Lemma \ref{barrierembedding} below and its proof appear in \cite{Todorcevicramseymethods} in a slightly different way, we add here a proof for completeness. 

\begin{nota}
    Given to families $\cB,\cC\subseteq[\w]^{<\w}$ we denote by $\cB\sqsubseteq\cC$ if
    \begin{itemize}
        \item $\forall s\in \cB\ \exists t\in\cC\ (s\sqsubseteq t)$,
        \item $\forall t\in \cC\ \exists s\in\cB\ (s\sqsubseteq t).$
    \end{itemize}
\end{nota}

\begin{lem}\label{barrierembedding}
    Given two barriers $\cB$ and $\cC$ on a countable set $N\in[\w]^\w$, there is an infinite set $M\in[N]^\w$ such that either $\cB|M\sqsubseteq\cC|M$ or $\cC|M\sqsubseteq\cB|M$. Moreover, if $\rho(\cB)<\rho(\cC)$ and $\cC$ is uniform, then $\cB|M\sqsubseteq\cC|M$ necessarily holds.
\end{lem}

\begin{proof}
    Define $\cB_0=\{b\in\cB:\exists c\in\cC(b\subseteq c)\}$. Since $\cB=\cB_0\cup(\cB\setminus\cB_0)$ we can find, by Nash-Williams theorem, an infinite set $M\in[N]^\w$ such that either, $(\cB\setminus\cB_0)|M=\emptyset$ or $\cB_0|M=\emptyset$. If $(\cB\setminus\cB_0)|M=\emptyset$ then we are done as this implies that $\cB|M\sqsubseteq \cC|M$. To see this notice that for every $b\in\cB|M$ ( and hence $b\in\cB_0$), there is $c\in\cC$ such that $b\subseteq c$, but if we define $X=b\cup (M\setminus\max(b))$, there is also $c'\in\cC$ such that $c'\sqsubseteq X$. As $b\subseteq c\in\cC$, it happens that $b\sqsubseteq c'$ and $c'\in\cC|M$. For $c\in\cC|M$, we can also find $b\in\cB|M$ such that $b\sqsubseteq c\cup(M\setminus(\max(c)+1))$. The case where $c\sqsubset b$ is not possible since $b\in\cB_0$ and thus there is another $c'\in\cC|M$ such that $b\subseteq c'$.

    Assume otherwise $\cB_0|M=\emptyset$ and define $\cC_0=\{c\in\cC|M:\exists b\in\cB|M(c\subseteq b)\}$.
    Again by Nash-Williams theorem, we can find $M'\in[M]^\w$ such that either $\cC_0|M'=\emptyset$ or $(\cC\setminus\cC_0)|M'=\emptyset$. The case $\cC_0|M'=\emptyset$ is impossible since we can find $b\in\cB|M'\subseteq\cB|M$ and $c\in\cC|M'$ such that $b,c\sqsubseteq M'$ and then either $b\subseteq c$ or $c\subseteq b$ contradicting the choices of $M$ and $M'$.
    Thus $(\cC\setminus\cC_0)|M'=\emptyset$ and as in the previous case we get that $\cC|M'\sqsubseteq \cB|M'$. 
    
    For the last assertion of the statement notice that if $\cC|M\sqsubseteq\cB|M$, then $T(\cC|M)\subseteq T(\cB|M)$. In particular,
    \begin{equation*}
    \rho(\cC)=\rho(\cC|M)=\rho_{T(\cC|M)}(\emptyset)\leq\rho_{T(\cB|M)}(\emptyset)=\rho(\cB|M)\leq\rho(\cB). \qedhere\end{equation*}
    \end{proof}

The previous lemma tells us that we can define a preorder on the family of all barriers on $\omega$ by defining $\cB\leq \cC$ if there is an infinite set $M$ such that $\cB|M\sqsubseteq\cC|M$. One may be tempted to say that the partial order defined by identifying $\cB$ and $\cC$ if $\cB\leq \cC$ and $\cC\leq\cB$, collapse to $\omega_1$ when restricted to the family of uniform barriers, that is, $\cB\leq\cC$ if and only if $\rho(\cB)\leq\rho(\cC)$. 
However, the unpleasant fact about this is that for two barriers $\cB$ and $\cC$ of the same rank, it is not always true that $\cB\leq\cC$ and $\cC\leq\cB$. 
For example, if $\cB$ is the Schreier barrier and $\cC$ is any uniform barrier defined such that $\rho_{T(\cC)}(\{n\})=f(n)$ for an increasing function $f\in\w^\w$ that is strictly bigger that the identity, then $\cC\nleq\cB$. Notice that any uniform barrier $\cB$ of rank $\w$ is completely determined by a function $f\in\w^\w$ that encodes the ranks of the first level on $T(\cB)$. A more complex but naturally defined coding of the ranks of the successors of any element $s\in T(\cB)$ also determines completely the structure of any barrier $\cB$.

We will see now that this function is the only obstruction and that, if one is willing to compress some finite intervals of $\omega$ into points, we can define a weaker relation between barriers that depends only on their ranks and extends $\leq$. This will allow us to prove, in some cases, that a property of barriers holds for all barriers of a given rank $\alpha$ (in particular a space being ``$\cB$-sequentially compact'') if and only if it holds for some single uniform barrier of the same rank (e.g., Corollary \ref{alphasequentiallyequivalence}).

\begin{dfn}
    Given two barriers $\cB$ and $\cC$ on a countable set $M$, we write $\cC\preceq\cB$ if there is a finite-to-one, non decreasing function $f\in\w^\w$ such that for every infinite subset $M'\in[M]^\w$ there exists $N\in[M']^\w$ so that $f\rest N$ is one-to-one and
    \begin{itemize}
        \item $\forall b\in(\cB|N)\ \exists c\in\cC\ (c\sqsubseteq f[b])$.
    \end{itemize}
\end{dfn}

Notice that if $\cB$ is uniform and $\rho(\cC)<\rho(\cB)$, then $\cC\preceq\cB$ with $f$ being the identity map, since $\cC|M\leq\cB|M$ for every restriction to an infinite set $M$. The condition that $f\rest N$ is one-to-one is superfluous as we can always shrink $N$ in order to get this property, but it will be convenient to add this to the definition to avoid saying it explicitly each time we use the preorder $\preceq$.\\
Having said this, it is clear that the relation $\preceq$ is transitive, that is, if $\cC\preceq\cB$ and $\cD\preceq\cC$ is witnessed by $f$ and $g$ respectively, then $\cD\preceq\cB$ and it is witnessed by $g\circ f$.

\begin{propo}\label{double arrow if uniformity and greater rank}
    If $\cB$ and $\cC$ are two barriers in $\omega$, $\cB$ is uniform and $\rho(\cB)\geq\rho(\cC)$, then $\cC\preceq\cB$.
\end{propo}

\begin{proof}

    Let us alternately define two sequences $\{a_i:i\in\omega\}$ and $\{b_i:i\in\omega\}$ such that $a_i<b_i<a_{i+1}$ for every $i\in\omega$.
    We will also denote by $\{k_i:i\in\omega\}$ the increasing sequence defined by these two sequences, that is, $k_{2i}=a_i$ and $k_{2i+1}=b_i$. We will also simultaneously define $f\in\w^\w$ which is completely determined by the sequence $\{k_i:i\in\w\}$ as follows: $f\rest[0,k_0)=0$ and $f(n)=i$ if and only if $n\in[k_i,k_{i+1})$. We shall show that $f$ satisfies the definition of $\cC\preceq\cB$. 
    To save notation, let us write $\rho_\cB$ and $\rho_\cC$ instead of $\rho_{T(\cB)}$ and $\rho_{T(\cC)}$ respectively. Also, given $s\in[\w]^\w\setminus T(\cC)$ we will define $\rho_\cC(s)=-1$ so that we can talk about the rank of $s$ even if $s$ is not in $T(\cC)$.
    
    \textbf{Step 0:}

    Since $\rho_\cB(\emptyset)\geq\rho_\cC(\emptyset)$, there are $a_0<b_0<\w$ such that $\rho_\cB(\{a_0\})\geq\rho_\cC(\{0\})$ and $\rho_\cB(\{b_0\})\geq\rho_\cC(\{1\})$. 
    
    Now assume we have defined $\{a_i:i\leq n\}$ and $\{b_i:i\leq n\}$.

    \textbf{Step a:}
    
\textit{\underline{Inductive hypothesis at $n$:}}

Let us call 
$$S_a(n)=\Big\{s\in T(\cB):s\subseteq\bigcup_{i< n}[a_i,b_i)\ \land\ \forall i< n(|s\cap [a_i,b_i)|\leq 1)\Big\}.$$

We will ensure that the following inductive hypothesis will hold throughout the construction: Given $s\in S_a(n)$

%\subseteq\bigcup_{i<n}[a_i,b_i)$ is such that $|s\cap[a_i,b_i)|\leq 1$ for every $i<n$, then

\begin{itemize}
    \item[$\phi_a(n,s)$] $\equiv\rho_\cB(s)\geq\rho_\cC(f[s])$
    \item[$\psi_a(n,s)$] $\equiv$ if $j\geq a_n$, then $\rho_\cB(s\cup\{j\})\geq\rho_\cC(f[s\cup\{2n\}])$. 
\end{itemize}

We will also write $\phi_a(n)$ as a shorthand for $\forall s\in S_a(n)\phi_a(n,s)$ and $\psi_a(n)$ as a shorthand for $\forall s\in S_a(n)\psi_a(n,s)$.
It is worth to point out that $S_a(n)$, %$L_a(n)$,
$\phi_a(n,s)$ and $\psi_a(n,s)$ do not mention $b_n$. Note that if $j\geq a_0$ we have that $\rho_\cB(\{j\})\geq\rho_\cB(\{a_0\})\geq\rho_\cC(\{0\})$ as $\cB$ being uniform. Thus $\psi_a(0)$ is satisfied and $\phi_a(0)$ is vacuously satisfies too (here $S_a(0)=\{\emptyset\}$).

\textit{\underline{Construction of $a_{n+1}$:}}

    We need to take care of all $s\in S_a(n+1)$ (note that this set is well defined since it only depends on $\{a_i,b_i:i\leq n\}$).
    
    First consider $s\in S_a(n)$. By $\phi_a(n,s)$ we have that $\rho_\cB(s)\geq\rho_\cC(f[s])$ and thus we can find $m_s^a$ such that 
    \begin{itemize}
        \item[$(\star^a_1)$] $\rho_\cB(s\cup \{m\})\geq \rho_\cC(f[s]\cup \{2n+2\})$ for every $m\geq m_s^a$. 
    \end{itemize}

    Otherwise $s\in S_a(n+1)\setminus S_a(n)$ and we can write $s'=s\setminus\{j\}$ where $j=\max(s)\geq a_n$ and $s'\in S_a(n)$. Then we know that $\rho_\cB(s'\cup\{j\})\geq\rho_\cC(f[s']\cup\{2n\})$ by $\psi_a(n,s)$. Since $j\in[a_n,b_n)$, we have that $f(j)=2n$. Thus 

    \begin{itemize}
        \item[$(\ast^a)$]\label{astlem3.5}     $\rho_\cB(s)=\rho_\cB(s'\cup\{ j\})\geq\rho_\cC(f[s']\cup \{f(j)\})=\rho_\cC(f[s'\cup\{j\}])=\rho_\cC(f[s])$
    \end{itemize}

    Then there is $m_s^a\in\w$ such that $\rho_\cB(s\cup \{m_s^a\})\geq \rho_\cC(f[s]\cup \{2n+2\})$. Since $\cB$ is uniform, 
    \begin{itemize}
        \item[$(\star^a_2)$] $\rho_\cB(s\cup \{m\})\geq \rho_\cC(f[s]\cup \{2n+2\})$ for every $m\geq m_s^a$. 
    \end{itemize}
    
    Finally define
    $$a_{n+1}=\max\{m_s^a\mid s\in S_a(n+1) \}.$$ 

    \textit{\underline{Inductive hypothesis for $n+1$:}}    
    
    We want to see that $\phi_a(n+1)$ and $\psi_a(n+1)$ hold. Let $s\in S_a(n+1)$. If $s\in S_a(n)$ then $\phi_a(n+1,s)$ follows from $\phi_a(n,s)$. In consequence we can assume that $s\in S_a(n+1)\setminus S_a(n)$ and let $s'=s\setminus\{j\}$ where $j=\max(s)$. Hence by
    $(\ast^a)$, we have that $\rho_{\cB}(s)\geq \rho_{\cC}(f[s])$ which means that $\phi_n(n+1,s)$ holds and consequently $\phi_a(n+1)$ does.
    
    Now let us see that $\psi_a(n+1)$ is also true. Let $s\in S_a(n+1)$ and $j\geq a_{n+1}$. If $s\in S_a(n)$ then $m_a^s$ satisfies $(\star_1^a)$ and if $s\in S_a(n+1)\setminus S_a(n)$ we have that $m_a^s$ satisfies $(\star_2^a)$, either case, as $a_n+1\geq m_a^s$ we have that $\psi_a(n+1,s)$ holds and so does $\psi_a(n+1)$.

\textit{\underline{Inductive hypothesis at $n$:}}

The construction is dual to that of part a. Let us start by calling 
$$S_b(n)=\Big\{s\in T(\cB):s\subseteq\bigcup_{i< n}[b_i,a_{i+1})\ \land\ \forall i<n(|s\cap [b_i,a_{i+1})|\leq 1)\Big\}.$$

The corresponding inductive formulas for $s\in S_b(n)$ are now:

\begin{itemize}
    \item[$\phi_b(n,s)$] $\equiv\rho_\cB(s)\geq\rho_\cC(f[s])$
    \item[$\psi_b(n,s)$] $\equiv$ if $j\geq b_n$, then $\rho_\cB(s\cup\{j\})\geq\rho_\cC(f[s\cup\{2n+1\}])$. 
\end{itemize}

We will write again $\phi_b(n)$ and $\psi_b(n)$ as shorthands for $\forall s\in S_b(n)\phi_a(n,s)$ and $\forall s\in S_b(n)\psi_a(n,s)$ respectively. In the particular case of $\phi_b$, we have that it is the same formula as $\phi_a$ with different domain. Since $a_{n+1}$ has already been defined, the interval $[b_n,a_{n+1})$ considered in $S_b(n+1)$ makes completely sense.

\textit{\underline{Construction of $b_{n+1}$:}}

    Given $s\in S_b(n)$, it follows from $\phi_b(n,s)$ that $\rho_\cB(s)\geq\rho_\cC(f[s])$ and thus we can find $m_s^b$ such that 
    \begin{itemize}
        \item[$(\star^b_1)$] $\rho_\cB(s\cup \{m\})\geq \rho_\cC(f[s]\cup \{2n+3\})$ for every $m\geq m_s^b$. 
    \end{itemize}

    Otherwise $s\in S_b(n+1)\setminus S_b(n)$ and we can write $s'=s\setminus\{j\}$ where $j=\max(s)\geq b_n$ and $s'\in S_b(n)$. Thus $\rho_\cB(s'\cup\{j\})\geq\rho_\cC(f[s']\cup\{2n+1\})$ by $\psi_b(n,s)$. Since $j\in[b_n,a_{n+1})$, we have that $f(j)=2n+1$. Then 

    \begin{itemize}
        \item[$(\ast^b)$]\label{astlem3.5b}     $\rho_\cB(s)=\rho_\cB(s'\cup\{ j\})\geq\rho_\cC(f[s']\cup \{f(j)\})=\rho_\cC(f[s'\cup\{j\}])=\rho_\cC(f[s])$
    \end{itemize}

    We can again find $m_s^b\in\w$ such that $\rho_\cB(s\cup \{m_s^b\})\geq \rho_\cC(f[s]\cup \{2n+3\})$ and moreover 
    \begin{itemize}
        \item[$(\star^b_2)$] $\rho_\cB(s\cup \{m\})\geq \rho_\cC(f[s]\cup \{2n+3\})$ for every $m\geq m_s^b$. 
    \end{itemize}
    
    We then define
    $$b_{n+1}=\max\{m_s^b\mid s\in S_b(n+1) \}.$$ 

    \textit{\underline{Inductive hypothesis for $n+1$:}}    
    
    If $s\in S_b(n)$ then $\phi_b(n+1,s)$ follows from $\phi_b(n,s)$ and otherwise $s\in S_b(n+1)\setminus S_b(n)$ can be written as $s'=s\setminus\{j\}$ where $j=\max(s)$. By
    $(\ast^b)$, we have that $\rho_{\cB}(s)\geq \rho_{\cC}(f[s])$ and then $\phi_b(n+1)$ holds.
    
    To see that $\psi_b(n+1)$ is true, let $s\in S_b(n+1)$ and $j\geq b_{n+1}$. Then $b_{n+1}\geq m_s^b$ and either by $(\star_1^b)$ or $(\star_2^b)$ we have that $\phi_b(n+1,s)$ holds.
    
    This finishes the construction and it remains to prove that $f$ works.\\

Clearly $f$ is finite-to-one and non-decreasing. Given $M\in[\w]^\w$, we can find an infinite subset $N\in[M]^\w$ such that $|N\cap[k_i,k_{i+1}]|\leq1$ for every $i\in\omega$ and either $N\cap\bigcup_{i\in\omega}[k_{2i},k_{2i+1})=\emptyset$ or $N\cap\bigcup_{i\in\omega}[k_{2i+1},k_{2i+2})=\emptyset$. Without loss of generality assume that $N\subseteq\bigcup_{i\in\omega}[k_{2i},k_{2i+1})=\bigcup_{i\in\omega}[a_i,b_i)$.

Let now $b\in\cB|N$. It follows from the choice of $N$ that there is an $n\in\omega$ such that $b\in S_a(n)$. By $\phi_a(n,b)$ we know that $0=\rho_\cB(b)\geq\rho_\cC(f[b])$. Necessarily $\rho_\cC(f[b])\in\{0,-1\}$, but either case we can conclude that there is $c\in\cC$ such that $c\sqsubseteq f[b]$. This finishes the proof.
\end{proof}

The previous lemma is new and it is the last piece of the theory of barriers that we need in order to prove that the general notion of barrier sequential compactness extending $n$-sequentially compact spaces, depends only on the rank of the associated uniform barrier. In the case that $\cB$ has rank $\w$, we can drop the requirement that the barrier is uniform.

\begin{coro}\label{doublearrowrankomega}
    Let $\cB$ and $\cC$ be barriers on $\w$ such that $\rho(\cC)\leq\rho(\cB)=\w$, then $\cC\preceq\cB$.
\end{coro}

\begin{proof}
    By the fact that the preorder $\preceq$ is transitive and the previous Proposition, it suffices to show that $\cS\preceq\cB$ where $\cS$ is the Schreier barrier.

    Since $\cB$ has rank $\w$, it follows from Lemma \ref{homogeneityofbarriers} that for every $i\in\w$, we can find $m_i$ such that $\rho_\cB(m)\geq i+1$ for every $m\geq m_i$. Let $f\in\w^\w$ be given by $f(n)=i$ if and only if $n\in[m_i,m_{i+1})$ (and we define $f\rest m_0$ arbitrarily). We claim that $f$ witnesses that $\cS\preceq\cB$.

    Indeed, given $M\in[\w]^\w$, we can find $N_0\in[M]^\w$ such that $N_0\cap[0,m_0)=\emptyset$ and $|N_0\cap[m_i,m_{i+1})|\leq 1$ for every $i\in\w$. 
    Let $n_0=\min(N_0)$ and for every $i\in\w$ define $n_{i+1}$ as follows: for each $j<i+1$ we can apply Corollary \ref{uniquenessoffiniterankbarriers} and find $k_j$ such that $\cB[n_j]=[\w\setminus k_j]^r$ for some $r\geq l$ where $l$ is the unique natural number such that $n_j\in[m_l,m_{l+1})$, pick $n_{i+1}>k_j$ for every $j<i+1$. Thus $N=\{n_i:i\in\w\}\in[M]^\w$.

    Fix $b\in(\cB|N)$ and let $n_i=\min(b)$. It follows from the definitions that $b\setminus\{n_i\}\in[\w\setminus k_j]^r$ for some $r\geq l$ where $n_i\in[m_l,m_{l+1})$. Then $|f[b]|\geq l+1$ as $f$ is one to one on $N$ but $f(n_i)=l$ and $\cS(l)=\{s\in[\w]^{<\w}:\min(s)=l\land|s|=l+1\}$. Thus we can conclude that $c\sqsubseteq f[b]$ defined as the first $l+1$ elements of $f[b]$ is an element of $\cS$.
\end{proof}

\section{Infinite dimensional sequential compactness}\label{sectionpositiveresults}

We now define the notion of $\cB$-sequential compactness for barriers $\cB$ and prove necessary results to formulate the natural notion of $\alpha$-sequential compactness, where $\alpha\in \omega_1$ corresponds to the rank of the barrier.  All the results of this section are natural generalizations of results 
 presented  in \cite{kubistopologicalramsey} for finite $n$.

%We will construct the counterexamples that show that the classes of $\alpha$-sequentially compact spaces and related spaces do not coincide (at least consistently, in some cases) giving analogous results to those presented in \cite{kubistopologicalramsey} and \cite{corralnRamsey}. 

\begin{dfn}
    Let $\cB$ be a barrier on $M\in[\w]^\w$, let $X$ be a topological space, $f:\cB\rightarrow X$ and $x\in X$. We say that $f$ \emph{converges} to $x$ if for all  $U\in\cN(x)$ there is $n\in\w$ such that $f[(\cB|(M\setminus n))]\subseteq U$.
\end{dfn}

The word ``converges'' might give rise to missunderstanding as $f$ has countable domain, however, the meaning of this word will be understood by the context. The function $f$ should be thought as a $\cB$-dimensional sequence or $\alpha$-dimensional sequence for $\alpha=\rho(\cB)$, rather than a classical 1-dimensional sequence indexed by the countable set $\cB$.

\begin{dfn}
    Let $\cB$ be a barrier on $\w$ and let $X$ be a topological space. We say that $X$ is \emph{$\cB$-sequentially compact} if for all $f:\cB\rightarrow X$ there is $M\in [\w]^\w$ such that $f\restriction (\cB|M)$ converges.
\end{dfn}

The natural stratification of barriers given by their ranks will give us a natural classification of $\cB$-sequentially compact spaces grouping them and associating them to a countable ordinal attending the following definition:

\begin{dfn}
    Let $\alpha<\omega_1$. We say that $X$ is $\alpha$\emph{-sequentially compact} if $X$ is $\cB$-sequentially compact for all barriers $\cB$ of rank $\alpha$.
\end{dfn}

The first application of the theory developed in Section \ref{barrierssection} and the main reason for the introduction of the order $\preceq$ is due to the next theorem:

\begin{thm}\label{double arrow implies B-seq implies C-seq}
    If $\cC\preceq\cB$ and $X$ is $\cB$-sequentially compact then $X$ is also $\cC$-sequentially compact.
\end{thm}

\begin{proof}
    Assume $X$ is $\cB$ sequentially compact and let $f\in\w^\w$ as in the definition of $\cC\preceq\cB$. Let $h:\cC\to X$ and define $\hat{h}:\cB\to X$ by $\hat{h}(b)=h(c)$ for any $c\in\cC$ that is $\sqsubseteq$-compatible with $f[b]$. Notice that since $\cC$ is a barrier, there is at least one of such $c$ and then $\hat{h}$ is well defined.

    Since $X$ is $\cB$-sequentially compact, there is $M\in[\w]^\w$ such that $\hat{h}\rest(\cB|M)$ converges to some $x\in X$. 
    Take $N\in[M]^\w$ such that
    $$\forall b\in(\cB|N)\ \exists c\in\cC\ (c\sqsubseteq f[b])$$
    given by the definition of $\cC\preceq\cB$ and define $M_0=f[N]$. 
    It turns out that $\hat{h}\rest(\cB|N)$ also converges to $x$.
    We shall prove that $h\rest(\cC|M_0)$ converges to the same $x\in X$.

    To see this fix an open neighborhood $U\in\cN(x)$. We can find $n\in\omega$ such that $\hat{h}(b)\in U$ whenever $b\in\cB|(N\setminus n)$. Take $c\in\cC|(M_0\setminus f(n))$. Thus $c=\{f(m_0),\ldots,f(m_j)\}$ for some set $a=\{m_i:i\leq j\}\subseteq N\setminus n$. Take any infinite set $N'\subseteq N$ such that $a\sqsubseteq N'$. We can find $b\in\cB|N$ such that $b\sqsubseteq N'$. It follows that $a\sqsubseteq b$ since otherwise
    % we would have that $b\sqsubset a$ and then, as $f\restriction N$ is injective, $f(b)\sqsubset f(a)=c$. On the other hand, as $b\in\cB|N$, there is $c'\in \cC$ such that $c\sqsubseteq f(b)$. Thus we have that $c'\sqsubseteq f(b)\sqsubset f(a)=c$, so $c'\sqsubset c$, which is impossible since both are in $\cC$. %This is the proof that $a\sqsubseteq b$, just in case
    $f[b]\sqsubset f[a]=c$ would contradict that $b\subseteq N$. Therefore $h(c)=\hat{h}(b)\in U$ as desired.
\end{proof}

\begin{coro}\label{alphasequentiallyequivalence}
    Let $X$ be a topological space and $\alpha\in\w_1$, the following are equivalent:
    \begin{enumerate}
        \item\label{equivalenceitem1} $X$ is $\alpha$-sequentially compact,
    
        \item\label{equivalenceitem2} $X$ is $\cB$-sequentially compact for every uniform barrier of rank $\alpha$,

        \item\label{equivalenceitem3} $X$ is $\cB$-sequentially compact for some uniform barrier $\cB$ of rank $\alpha$.
    \end{enumerate}
\end{coro}

\begin{proof}
    That (\ref{equivalenceitem3}) implies (\ref{equivalenceitem1}) follows directly from Proposition \ref{double arrow if uniformity and greater rank} and Theorem \ref{double arrow implies B-seq implies C-seq}. The other implications are clear from the definitions.
\end{proof}

From Lemma \ref{double arrow if uniformity and greater rank} and Theorem \ref{double arrow implies B-seq implies C-seq} we can easily get the following result.

\begin{coro}\label{downwardsimplicationseuqntialcompactness}
    If $X$ is $\alpha$-sequentially compact and $\beta<\alpha$ then $X$ is $\beta$-sequentially compact. In particular, if $\alpha$ is infinite, $X$ is $n$-sequentially compact for every $n\in\omega$.
\end{coro}

The stronger result for barriers of rank $\w$ given in Corollary \ref{doublearrowrankomega} yields the following result:

\begin{coro}
    If $X$ is $\cB$-sequentially compact for some barrier of rank $\omega$ (not necessarily uniform), then $X$ is $\omega$-sequentially compact.
\end{coro}

In general, this is not true for every countable ordinal. A trivial example is a barrier with a single node in $T(\cB)$ of rank $\w$ (let us say $\{0\}$) and such that $\cB|(\omega\setminus1)$ is isomorphic to $\cS_2$. In this case, every infinite restriction has rank $\w$ but the barrier itself has rank $\w+1$. It is not hard to modify our constructions in Section \ref{counterexamplessection} to show that, consistently, there is a $\cB$ sequentially compact space that is not $(\w+1)$-sequentially compact. We do not know if the requirement that $\cB$ has an infinite restriction of rank $\alpha$ suffices in order to show that $\cB$-sequentially compact spaces are $\alpha$-sequentially compact. We believe that the answer is ``no'' but we conjecture that a weaker result is true.

\begin{con}\label{conjecture of terminal rank}
If $X$ is $\cB$-sequentially compact, then it is $\rho_u(\cB)$-sequentially compact.
\end{con}

The following easy fact is a partial answer to Conjecture \ref{conjecture of terminal rank}.

\begin{propo}
    Let $\cB$ be a barrier on $\w$ and $X$ a $\cB$-sequentially compact space. If there is $M\in[\w]^\w$ such that $\cB|M$ is uniform and either $X$ is not $(\cB|(\w\setminus M))$--sequentially compact or $\w\setminus M$ is finite, then $X$ is $\rho(\cB|M)$-sequentially compact. In particular $X$ is $\rho_u(\cB)$-sequentially compact.
\end{propo}

\begin{proof}
    By Corollary \ref{alphasequentiallyequivalence} it is enough to prove that $X$ is $(\cB|M)$-sequentially compact. 
    Let $F:=\w\setminus M$. Regardless $F$ is finite or $X$ is not $(\cB|F)$-\seqcompp, %we know that  %Let $f:\cB|M\to X$, and fix $\hat{f}:\cB\to X$ any extension of $f$. As $X$ is $\cB$-sequentially comopact, there is $N\in[\w]^\w$ such that $\hat{f}\restriction N$ is convergent, but then $f|\restriction (N\setminus F)$ converges to the same point.
    %Case 2) $X$ is not $\cB|(\w\setminus M)$-\seqcomp. Then 
    there exists a function $g:\cB|F\to X$ without infinite convergence subsequences (here convergence means convergence with respect to barriers).
    %\footnote{Note that if $\cB|F=\emptyset$ then $g$ is the empty function.}. 
    Now let $f:\cB|M\to X$. We want to prove that $f$ admits an infinite convergent subsequence $f\rest(\cB|M')$ for some $M'\in[\w]^\w$. 
    For this let $\hat{f}:\cB\to X$ be any common extension of both $f$ and $g$. As $X$ is $\cB$-\seqcompp, there is $N\in[\w]^\w$ such that $\hat{f}\restriction N$ is convergent. Now $N\cap F$ should be finite, since otherwise $\hat{f}\restriction (N\cap F)$ would also be an infinite convergent subsequence for $g$. This way $N':=N\setminus F$ is infinite and thus $\hat{f}\restriction N'=f\restriction N'$ is an infinite convergent subsequence for $f$.  
\end{proof}

\section{Some classes of $\omega_1$-sequentially compact spaces}\label{examples}

In \cite{kubistopologicalramsey}, it is proved that compact metric spaces are $n$-sequentially compact for every $n\in\w$. This result extends to any $\alpha<\w_1$. We will start this section by giving two classes of spaces, each containing the class of compact metric spaces, that sit at the top of the hierarchy of $\alpha$-sequentially compact spaces.

\begin{dfn}
    We say that $X$ is \emph{$\w_1$-sequentially compact}, if it is $\alpha$-sequentially compact for every $\alpha<\w_1$.
\end{dfn}

The first class of spaces that we want to show to be $\w_1$-sequentially compact, is the class of sequentially compact spaces of character less that $\b$. The fact that sequentially compact spaces of character less that $\b$ are $n$-sequentially compact for every $n\in\w$, was proved in \cite{kubistopologicalramsey}.

\begin{thm}\label{characterlessthanbimpliesBseq} Suppose $X$ is a sequentially compact space. If $\chi(X)<{\mathfrak b}$ then $X$ is $\w_1$-sequentially compact. %$\cB$-sequentially compact for every barrier $\cB$.
\end{thm}

\begin{proof}
We will show by induction that $X$ is $\alpha$-sequentially compact for every $\alpha<\omega_1$. If $\alpha=1$, then $X$ being sequentially compact is equivalent to $X$ being 1-sequentially compact as the only uniform barrier of rank $1$ is $[\omega]^1$.

We proceed to prove the inductive step at $\beta$, regardless of it being a successor ordinal, and hence of the form $\beta=\alpha+1$ or a limit ordinal which would then be the $\operatorname{sup}\alpha_n$, where $\alpha_n$ form an increasing sequence. Now let $f:\cB\to X$ be any function where $\cB$ is uniform of rank $\beta$.
For all $n\in\omega$, $\cB[n]$ is either an $\alpha$-uniform barrier on $\omega\setminus (n+1)$ or an $\alpha_n$-uniform barrier if $\beta$ were a limit ordinal. 

For all $n\in\omega $ let $f_n:\cB[n]\to X$ be the function mapping $s\mapsto f(\{n\}\cup s)$. By our inductive hypothesis there is an infinite subset $N_0$ and a point $x_0\in X$ such that $f_0\restriction(\cB[0]|N_0)\rightarrow x_0$. %\carlos{I think is a little bit more precise the notation $f_0\restriction(\cB_{\{0\}}| N_0)\rightarrow x_0$ instead of $f_0[\cB_{\{0\}}| N_0]\rightarrow x_0$. Also at some point $\cB_{\{n\}}$ has become simple $\cB_n$, so im adding the $\{$ and $\}$}. 
Let $n_0=0$ and $n_1=\min N_0\setminus\{0\}$.

Note that for any infinite subset $A$ of $\omega$ and any $n\in\omega$, $\cB[n]|A$ is also uniform of rank $\alpha$ if $\beta=\alpha+1$ or of rank $\alpha_n$ if $\beta$ were a limit ordinal.  
Now suppose that we have defined a decreasing sequence of infinite subsets $N_0\supseteq N_1\supseteq\dots\supseteq N_{i-1}$, and an increasing sequence $n_0<n_1<\dots<n_i$ in addition to a subset $\{x_0,\dots,x_{i-1}\}$ of $X$. 
Consider the function $f_{n_{i}}:(\cB[n_i]| N_{i-1})\rightarrow X$, applying our inductive hypothesis yields an infinite subset $N_{i}\subseteq N_{i-1}$ such that $f_{n_i}\restriction(\cB[n_i]|N_{i})\rightarrow x_i$, for some point $x_i\in X$. We also define $n_{i+1}=\min(N_i\setminus\{n_0,\ldots,n_i\})$. Thereby, we obtain a decreasing sequence of infinite subsets $N_0\supseteq N_1\supseteq N_2\supseteq\dots$, a subset $\{x_i:i\in\omega\}$ of the space $X$  and a strictly increasing sequence $\langle n_i:i\in\omega\rangle$ such that $f_{n_i}\restriction(\cB[n_i]|N_i)\rightarrow x_i$ where $n_{i+1}=\min N_i\setminus\{n_0,\ldots,n_i\}$ for all $i\in\omega$.

Since $\{x_i:i\in\omega\}\subseteq X$ and $X$ is sequentially compact, there is a convergent subsequence, $Y=\{x_{i_j}:j\in\omega\}$, with limit point $x.$ Re-numerating the indices, we can view $Y=\{x_i:i\in\omega\}$ as the convergent subsequence. Also, let $N=\{n_i:i\in\omega\}$

Take any open set $U$ in $X$ containing the point $x$. As $x_i\to x$ there is an integer $m_U\in\omega$ such that for all $i>m_U$ we have that $x_i\in U.$ Also, for all such $i$ there is another integer $\phi_U(n_i)$ such that 
\[
f_{n_i}[\cB[n_i]| (N_i\setminus\phi_U(n_i))]\subseteq U.\tag{$\ast$}
\] \label{star}

For any $j\leq m_U$ and any $j\notin\{n_i:i\in\omega\}$ let $\phi_U(j)=0$. Thus, any open subset $U$ of $X$ containing $x$, induces a function $\phi_U:\omega\rightarrow\omega$.

Let $\eta (x)$ be a local neighborhood base of $x$ of minimum size. Since $\chi(X)<\mathfrak b$ the set of functions $\{\phi_U:U\in\eta(x)\}$ is bounded by some increasing function $\psi:\omega\rightarrow\omega$, that is $\phi_U<^*\psi$ for every $U\in\eta(x)$.

We can find an increasing subsequence $\langle m_i:i\in\omega\rangle$ of $\langle n_i:i\in\omega\rangle$ such that $m_{i+1}>\psi(m_i)$ for every $i\in\omega$. We claim that $f\rest(\cB|M)\to x$ where $M=\{m_i:i\in\w\}$.  

Fix an open set $U\in\eta(x)$ and let $k\in\omega$ such that $x_j\in U$ and $\psi(j)>\phi_U(j)$ for every $j\geq k$. Let $s\in\cB|(M\setminus \psi(k))$ and let us write $s=\{m_{k_0},\ldots,m_{k_n}\}$ and $s'=s\setminus\{m_{k_0}\}$. Note that for every $i>0$, we have that $$m_{k_i}>\psi(m_{k_0})>\phi(m_{k_0}),$$
as $m_{k_0}\geq\psi(k)\geq k$. In particular, $s'\cap\phi(m_{k_0})=\emptyset$.

Now let $j\in\omega$ such that $m_{k_0}=n_j$. Note that if $j'>j$, then $n_{j'}\in N_j$. As any $m_{k_i}>m_{k_0}$ for $i>0$ and this implies that $m_{k_i}=n_{j'}$ for some $j'>j$, we can conclude that $s'\subseteq N_j$.

Combining the previous arguments we get that $s'\subseteq N_j\setminus\phi_U(n_j)$ and we already know that $s=\{n_j\}\cup s'\in\cB$, which implies that $s'\in\cB[n_j]$. Therefore $f(s)=f_{n_j}(s')\in U$ by ($\ast$) as desired.
\end{proof}

The second class of spaces that we shall show are $\w_1$-sequentially compact is the class of compact bisequential spaces.
Recall that $\cU\subseteq\cP(X)$ clusters at $x$ (often written as $x\in\overline\cU$) if $x\in\overline U$ for every $U\in\cU$. Of course if $\cU$ is an ultrafilter, then $\cU$ clusters at $x$ is equivalent to $\cU$ converging to $x$.

\begin{dfn}\cite{michael1972quintuple}
    Let $X$ be a topological space and $x\in X$. We say that $X$ is \emph{bisequential at $x$}, if for every ultrafilter $\cU$ that converges to $x$, there is a countable subfamily $\{A_n: n\in\w\}\subseteq\cU$ such that for every $W\in\cN(x)$ there is $n\in\omega$ such that $A_n\subseteq W$. We say that $X$ is \emph{bisequential} if it bisequential at $x$ for every $x\in X$.
\end{dfn}

We say that $X$ is {\em countably bisequential} if for every ultrafilter on a countable subset of $X$ that clusters at $x\in X$, there is a countable subfamily of the ultrafilter converging to $x$.

Given a tree $T\subseteq\w^{<\w}$ and $s\in T$, recall that the set of successors of $s$ in $T$ is $succ_T(s)=\{n\in\w:s\conc n\in T\}$. Given $\cU\subseteq\cP(\w)$, we say that $T$ is $\cU$\emph{-branching} if $succ_T(s)\in\cU$ for every $s\in T$. As usual, we will omit the subindex in $succ_T$ when there is no risk of confusion. If $\cU$ is an ultrafilter and $\cB$ is a barrier on $\w$, we define an ultrafilter $\cU^\cB\subseteq \cP(\cB)$ by declaring $Y\in\cU^\cB$ if and only if, there is an $\cU$-branching tree $S\subseteq T(\cB)$ such that $Y=S\cap \cB$ (here we are identifying finite subsets with finite increasing sequences on $\w$).%\cesar{Should we mention here the other possible (but equivalent) definitions for $\cU^\cB$? If so, we need to look for references for them and the proofs of their equivalences.}

\begin{thm}\label{bisequentialimpliesramsey}
Every compact countably bisequential space is $\w_1$-sequentially compact.
\end{thm}

\begin{proof}
    Let $X$ be a compact countably bisequential space, let $\cB$ be a barrier and let $f:\cB\to X$. Fix an ultrafilter $\cU\subseteq\cP(\w)$ and define an ultrafilter $\cV\subseteq \cP(f[\cB])\subseteq\cP(X)$ by $V\in\cV$ if and only if $f^{-1}(V)\in\cU^\cB$.

    As $X$ is compact, we can find $x\in X$ such that $x\in\overline\cV$ and by bisequentiality, there exists a decreasing family $\{V_n:n\in\w\}\subseteq\cV$ such that for any open neighborhood $W\in\cN(x)$, there is $n\in\w$ such that $V_n\subseteq W$. 
    Define $U_n=f^{-1}(V_n)$ for every $n\in\w$ and let $T_n=\{s\in T(\cB):\exists b\in U_n(s\subseteq b)\}$. Then each $T_n$ is an $\cU$-branching tree. 
    
    We now recursively define an infinite set $M=\{m_i:i\in\omega\}$. Let $m_0=\min(succ_{T_0}(\emptyset))$. If we have already defined $\{m_i:i<n\}$, define $m_n$ such that for every $s\subseteq\{m_i:i<n\}$ with $s\in T(\cB)\setminus\cB$, we have that $s\conc m_n\in T_n$. The choice of $m_n$ is possible since every relevant $s$ is an element of $T_0$ and the trees $T_n$ are $\cU$-branching. Moreover, if $s\subseteq M$, $s\in T(\cB)\setminus\cB$ and $\min(s)\geq m_n$, then $s\in T_k$ for every $k\geq n$.

    It remains to show that $f\rest(\cB|M)$ converges to $x$. Fix an open set $W\in\cN(x)$%\carlos{i changed the notation of $\eta(x)$ for $\cN(x)$} 
    and find $n\in\omega$ such that $V_n\subseteq W$. Given $s\in\cB|(M\setminus m_n)$, let $m_k=\max(s)$ and $s'=s\setminus m_k$. 
    By the choice of $m_k$, we have that $s= {s'}\conc m_k \in T_n\cap \cB=U_n$
    as $k\geq n$, thus $f(s)\in V_n\subseteq W$. 
    Therefore $f[\cB|(M\setminus m_n)]\subseteq W$.
    \end{proof}

Recall that a function is of Baire-class-1 if it is the pointwise limit of a sequence of continuous functions. A \emph{Rosenthal compactum} is a compact subset $K\subseteq B_1(X)$, where $B_1(X)$ denotes the set of all Baire-class-1 functions from $X$ to $\bR$ endowed with the pointwise topology, for some Polish space $X$. Roman Pol proved that Rosenthal compacta are countably bisequential \cite{Pol1989note}. 
It can also be deduced from arguments of G. Debs in \cite{debseffective} (see \cite{todorcevicramseyspaces} and Lemma 6 in \cite{todorcevicbairefirstclass}).
Therefore we have 
\begin{coro}
    Every Rosenthal compact space is $\w_1$-sequentially compact.
\end{coro}

A space $X$ is \emph{angelic} if relatively countably compact subsets of $X$ are relatively compact and for every relatively compact subset $A\subseteq X$ and every $x\in\overline{A}$ there is a sequence $\{x_n:n\in\omega\}\subseteq A$ that converges to $x$. It is a result of  J. Bourgain, D. H. Fremlin and M. Talagrand \cite{bourgainfremlintalagrand} that Rosenthal compacta are angelic spaces.

It was shown by Knaust in \cite{knaustarray} that Rosenthal compacta have the Ramsey property. We recall that a space $X$ has the Ramsey property if for every function $f:[\w]^2\to X$ such that $\lim_{i\to\infty}\lim_{j\to\infty}f(\{i,j\})=x$, there is $M\in[\w]^\w$ such that $f\rest[M]^2\to x$. As any 2-sequentially compact space has the Ramsey property, the above corollary is a strengthening of Knaust's result.

The lists of angelic spaces with the Ramsey property was expanded in \cite{knaustangelic}, however the question of whether every angelic space has the Ramsey property, was left open. We show now that this is not the case by pointing out that the example of a sequentially compact space that is not 2-sequentially compact considered in \cite{corralnRamsey} is such a counterexample.

\begin{thm}
    There is an angelic space without the Ramsey property.
\end{thm}

\begin{proof}
    Note that a compact space is angelic if and only if it is Fr\'echet. In \cite{corralnRamsey} an almost disjoint family $\cA\subseteq[\w\times\w]^\w$ such that $\cF(\cA)$ is Fr\'echet, sequentially compact but fails to be 2-sequentially compact is constructed. Hence $\cF(\cA)$ is an angelic space.

    In order to prove that $\cF(\cA)$ is not 2-sequentially compact, $\cA$ was constructed so that $A_n:=\{(n,m):m\in\omega\}\in\cA$ for every $n\in\w$ and so that $G:[\w]^2\to\w\times\w$ does not converge to any point in $\cF(\cA)$. To avoid confusion, let's denote by $\ast$ the point at infinity on $\cF(\cA)$.
    As $\lim_{m\to\infty}(n,m)=A_n$ and any infinite subset of $\cA$ converges to $\ast$ in $\cF(\cA)$, we get that $\lim_{n\to\infty}\lim_{m\to\infty}=\ast$. Then moreover we have that $\cF(\cA)$ does not have the Ramsey property. 
\end{proof}

\vspace{3mm}
The notion of $n$-sequentially compact space in \cite{kubistopologicalramsey} was motivated by the 2-dimensional version introduced in \cite{ramseymetric}, where the main application was to show that compact metric semigroups have idempotents naturally defined as the limit of a 2-dimensional sequence.
From Theorem \ref{bisequentialimpliesramsey} and the results in \cite{ramseymetric}, we can conclude the following:

\begin{coro}
Every compact countably bisequential semigroup $K$ has an idempotent naturally
representable as the Ramsey limit of a restriction of any given $f:[\w]^2\to K$.  
\end{coro}

\vspace{2mm}
The examples of $n$-sequentially compact spaces that are not $n+1$-sequentially compact that were constructed in \cite{kubistopologicalramsey} and \cite{corralnRamsey} were all of the form ${\mathcal F}({\mathcal A})$, i.e., the one point compactifications of $\Psi$-spaces over the almost disjoint family ${\mathcal A}$. An important class of almost disjoint families are those whose $\Psi$-space can be continuously embedded into ${\mathbb R}$. 
An almost disjoint family $\cA$ on $\w$ is $\bR$\emph{-embeddable} if there is a one-to-one function $f:\w\to\bQ$ that can be continuously extended to a function $\widehat{f}:\Psi(\cA)\to\bR$ (see \cite{GuzmanOnRembeddability} and \cite{hrusakqsetsandnormality}). 
%Given $X\subseteq2^\w$, define $\cA_X=\{\{x\rest n:n\in\omega\}:x\in X\}$. We call such a family an \emph{almost disjoint family of branches}. It is easy to see that $\cA$ is $\bR$-embeddable if $\Psi(\cA)$ is homeomorphic to $\Psi(\cA_X)$ for some $X\subseteq2^\w$.
We now note that the $\cR$ embeddable families all give rise to $\omega_1$-sequentially compact spaces ${\mathcal F}({\mathcal A})$. As a $\Psi$-space is first countable and hence bisequential, when we say that an almost disjoint family is bisequential, we will mean that $\cF(\cA)$ is.
%(we recall that $\cF(\cA)$ is the one-point compactification of $\Psi(\cA)$).
By Proposition 3.2 in \cite{gruenhageFUFin2} every $\bR$-embeddable family is bisequential, so we have:

\begin{coro}
    Every $\bR$-embedable almost disjoint family is $\w_1$-sequentially compact, i.e., if $\cA$ is $\bR$-embeddable then $\cF(\cA)$ is $\w_1$-sequentially compact.
\end{coro}

\vspace{3mm}

A compact space is said to be \emph{(uniform) Eberlein compact} if it embeds homeomorphically into a Banach (Hilbert) space with its weak topology. A compact space is \emph{Corson compact} if it homeomorphically embeds into a $\Sigma$-product, $\Sigma(\kappa)$, where
$$\Sigma(\kappa)=\{f\in\bR^\kappa:|supp(f)|\leq\w\}$$
and $supp(f)=\{\alpha\in\kappa:f(\alpha)\neq 0\}$. A nice way to obtain Eberlein and Corson compactums is via adequate families of subsets introduced by Talagrand \cite{Talagrand}.
A family $\cD\subseteq \cP(X)$ is \emph{adequate} if it is closed under subsets and $D\in\cD$ if and only if every finite subset of $D$ is an element of $\cD$. Every adequate family of subsets of $X$ defines a compact subset of $2^X$ by using characteristic functions. 
We call the set $K_\cD=\{1_D\in2^X:D\in\cD\}$ an \emph{adequate compact}, where $1_D$ is the characteristic function for $D$. 
It is clear that $K_\cD$ is Corson if and only if $\cD\subseteq[X]^{\leq\w}$.
A characterization of (uniform) Eberlein adequate compacts is due to A. G. Leiderman and G. A. Sokolov.

\begin{thm}\cite{leidermanadequatefamilies}
    Let $\cD\subseteq\cP(X)$ be an adequate family, then:
    \begin{enumerate}[label=(\roman*)]
        \item $K_\cD$ is Eberlein compact if and only if $X=\bigcup_{n\in\omega}X_n$ and $supp(f)\cap X_n$ is finite for every $f\in K_\cD$ and $n\in\w$. 
        \item $K_\cD$ is uniform Eberlein compact if and only if $X=\bigcup_{n\in\omega}X_n$ and there is a function $g\in\w^\w$ such that $|supp(f)\cap X_n|<g(n)$ for every $f\in K_\cD$ and $n\in\w$.
    \end{enumerate}
\end{thm}

In order to see that Eberlein compacta defined from adequate families are $\w_1$-sequentially compact we need the following definition.

\begin{dfn}\cite{arkhangelskiifunctionspaces}
    A space $X$ is \emph{strongly monolithic} if the weight of $\overline{A}$ does not exceed $|A|$ for every $A\subseteq X$.
\end{dfn}

We are only interested in the case where $|A|=\w$. Note that  (uniform) Eberlein compacta of the form $K_\cD$ for an adequate family $\cD$ and Corson compacta are strongly monolithic. Hence the next result will show that they are also $\w_1$-sequentially compact.

\begin{thm}\label{eberlein}
Every sequentially compact strongly monolithic space is $\w_1$-sequentially compact. In particular, Corson compacta and any Eberlein compacta defined from adequate families are $\w_1$-sequentially compact.    
\end{thm}

\begin{proof}
    As the closure of any countable set on a strongly monolithic space $X$ is first countable, we can apply Theorem \ref{characterlessthanbimpliesBseq} to $f[\cB]$ for every barrier $\cB$ in order to find a convergent subsequence. Hence $X$ is $\w_1$-sequentially compact.
    As Corson and Eberlein compacta are sequentially compact \cite{eberleinweakcomppactness}, the second part of the result follows.
\end{proof}

\section{Some examples delineating the classes of $\alpha$-sequentially compact spaces.}\label{counterexamplessection}

We will construct the counterexamples that show that the classes of $\alpha$-sequentially compact spaces and related spaces do not coincide (at least consistently, in some cases) giving analogous results to those presented in \cite{kubistopologicalramsey} and \cite{corralnRamsey}. 

We shall start this section by pointing out that neither of the theorems \ref{characterlessthanbimpliesBseq} or \ref{bisequentialimpliesramsey} supercede the other. On the one hand, if $X\subseteq 2^\w$ has size $\c$, then $\cF(\cA_X)$ is a compact bisequential space but its character is $|\cA_X|=\c\geq\b$ (so satisfies the hypotheses of \ref{characterlessthanbimpliesBseq} but not \ref{bisequentialimpliesramsey}). On the other hand, $\w_1$ is first countable, sequentially compact spaces but not compact (so satisfies \ref{bisequentialimpliesramsey} but not \ref{characterlessthanbimpliesBseq}). However, this leaves the question whether there is an $\omega_1$-sequentially compact, compact space that is not bisequential, or even one of character less than $\b$ (of course, assuming $\b>\omega_1$ since first countable spaces are bisequential).
%Note that, since $\b=\w_1$ is consistent, spaces of character less than $\b$ are consistently bisequential so we need to assume $\b>\omega_1$.

\vspace{3mm}

In order to find a compact and $\w_1$-sequentially compact space that is not bisequential, we will use the following space defined in \cite{Ciesla}: consider the set
$$\cD=\{A\in[\w_1\times\w_1]^{<\w}:A\textnormal{ is the graph of a partial strictly decreasing function}\}.$$
Clearly $\cD$ is an adequate family and then $\cZ=K_\cD$ is a subset of the Tychonoff product $2^{\w_1\times\w_1}$ by identifying elements of $\cD$ with their characteristics functions. Then $\cZ$ is an Eberlein compactum of weight $\omega_1$ given by an adequate family.

\begin{thm}\cite{Ciesla}
The space $\cZ$ is not bisequential.    
\end{thm}

With the previous result and Theorem \ref{eberlein} we get the following:

\begin{thm}
    There exists a compact and $\w_1$-sequentially compact space that is not bisequential.\QED
\end{thm}

If we moreover assume that $\b>\omega_1$, then the example has character $<\b$.

A more general version of $\cZ$ was considered by C. Agostini and J. Somaglia in \cite{agostinibisequentiality}. Given a tree $T$, we say that a finite partial function $f;T\to T$ is \emph{decreasing} if $s<t$ implies $f(s)>f(t)$. It is \emph{strongly decreasing}, if $f$ is decreasing and $dom(f)$ is a totally ordered subset of $T$. Hence $\cZ_{T}\subseteq2^{T\times T}$ is the space of strongly decreasing finite partial functions on $T$ and $\cY_T\subseteq2^{T\times T}$ is the space of finite decreasing functions on $T$.

\begin{thm}
    If $T$ is a tree with no uncountable antichains and it has an uncountable branch, then $\cZ_T$ and $\cY_T$ are compact and $\w_1$-sequentially compact spaces but fail to be bisequential.
\end{thm}

\begin{proof}
    Bisequentiality of $\cZ_T$ and $\cY_T$ was characterized in terms of combinatorial properties of $T$.
The space $\cZ_T$ is bisequential if and only if $T$ has size less than the first measurable cardinal and $T$ has no uncountable branches whilst $\cY_T$ is bisequential only in the case that $T$ is countable \cite{agostinibisequentiality}.
It is also easy to see that $\cZ_T$ is always Eberlein compact (defined from an adequate family) and that $\cY_T$ is Corson compact if and only if every antichain in $T$ is countable.
\end{proof}

%\begin{example}
%There is an $\w_1$-sequentially compact space that is not bisequential.    
%\end{example}

%\begin{thm}
%    The following are equivalent:
  %  \begin{enumerate}
 %       \item $\b=\w_1$
 %       \item every sequentially compact space of character less that $\b$ is bisequential
%    \end{enumerate}
%\end{thm}

%Notice that if $\b>\w_1=\chi(\cY_T)$, the previous example also shows that not every sequentially compact space of character less that $\b$ is bisequential. On the other hand, if $\b=\w_1$, then $\cY_T$ is an $\w_1$-sequentially compact space of character $\b$. 

%\begin{que}
%    Is every sequentially compact Fr\'echet space of character less than $\b$ bisequential?
%\end{que}

%\begin{que}\label{questioncountabletightnessandsequential}
%What if we replace Fr\'echet for countable tightness or if the space is sequential instead? 
%\end{que}

%Remind that in the presence of $\textsf{PFA}$ every compact space of countable tightness is sequential \cite{Baloghcountabletightness}, hence both modifications in Question \ref{questioncountabletightnessandsequential} are consistently equivalent for compact spaces.

It is worth noting that for an almost disjoint family $\cA$, if its Franklin space $\cF(\cA)$ is $2$-sequentially compact, then $\cA$ is  nowhere mad. Otherwise, if $\cA\rest X$ is mad with $X\in\cI(\cA)^+$, we can find a function $f:[\omega]^2\to X$ that does not converge in $\cF(\cA)$ as mad families are never $2$-sequentially compact (see \cite{kubistopologicalramsey}). Since an ad family is Fr\'echet if and only if it is nowhere mad, all counterexamples constructed from ad families, that are at least $2$-sequentially compact in \cite{kubistopologicalramsey} and \cite{corralnRamsey} are also Fr\'echet. It was also shown in \cite{corralnRamsey} that 2-sequentially compact spaces are $\alpha_3$.
This gives more examples of non-bisequential almost disjoint families that are Fr\'echet and $\alpha_3$. The existence of these kind of ad families is not known in \textsf{ZFC}. For more consistent examples and for the definition of $\alpha_3$ the reader may  consult \cite{corralfrechetlike}.\\

It was previously shown in \cite{kubistopologicalramsey} that the classes of $n$-sequentially compact spaces and $(n+1)$-sequentially compact spaces does not coincide (assuming \textsf{CH} for $n>1$). This result was later improved by the authors in \cite{corralnRamsey}, by showing that the same holds under weaker assumptions, e.g., $\b=\c$. 
We will show that we can also consistently differentiate between $\alpha$-sequentially compact and $\beta$-sequentially compact spaces for any $\alpha\not=\beta$, by stepping up the combinatorial analysis done in \cite{corralnRamsey} to barriers.\\

We will say that $T\subseteq\omega^{<\w}$ is a {\it Hechler tree}, if for every $s\in T$, the set of successors of $s$ is \emph{coinitial} on $\w$, i.e., there exists $k_s\in\omega$ such that $succ_T(s)=\{i\in\omega: s\conc i\in T\}=\omega\setminus k_s$. We can naturally associate a function $h_T:T\to\omega$ to each Hechler tree $T$ by defining $h_T(s)=k_s$. 

For $s,t\in\omega^{<\w}$, we say that $s\prec t$ if either $s\sqsubset t$ or $s=r\conc l$, $t=r\conc m$ for some $r\in\omega^{<\w}$ and $l<m$. Here $s\sqsubset t$ means that $s$ is a proper initial segment of $t$.
Enumerate $\omega^{<\w}=\{s_i: i\in\omega\}$ such that $s_i\prec s_j$ implies $i<j$. We fix this enumeration for the rest of the paper.

\begin{dfn}
    Let $f\in\w^\w$ and let $T\subseteq\w^{<\w}$ be a Hechler tree. We define:
    \begin{itemize}
        \item $T_f\subseteq\omega^{<\w}$ where $\emptyset\in T_f$ and $succ(s_i)=\omega\setminus f(i)$ for every $s_i\in T_f$.
        \item $f_T\in\omega^\omega$ where $f_T(n)=k$ if and only if $s_n\in T$ and $succ_T(s_n)=\omega\setminus k$ and $f_T(n)=0$ otherwise.
    \end{itemize}
\end{dfn}

Recall that $\textsf{FIN}\subseteq\cP(\w)$ is the ideal of finite sets of $\w$.

\begin{dfn}
    Given a barrier $\cB$, define the ideal $\textsf{FIN}^\cB$ as the set of all $X\subseteq[\w]^{<\w}$ such that
    $$\{n\in\omega:\{s\setminus\{n\}:s\in X\cap\cB(n)\}\notin\textsf{FIN}^{\cB[n]}\}\in\textsf{FIN}.$$
\end{dfn}

Notice that in the case of the barrier $\cB=[\w]^n$, the ideal $\textsf{FIN}^\cB$ coincides with the well known Fubini product $\textsf{FIN}^n$. In general, the Fubini product $\textsf{FIN}^\alpha$ for a countable limit ordinal $\alpha$ is not uniquely determined and depends on the choice of an increasing sequence $\{\alpha_n:n\in\w\}$ converging to $\alpha$, the corresponding previous choices for each $\alpha_n$ and so on. For a barrier $\cB$ of rank $\alpha$ such that $\alpha_n=\rho(\cB[n])$ forms an increasing sequence, the ideal $\textsf{FIN}^\cB$ also coincides with a Fubini product $\textsf{FIN}^\alpha$ where the sequence $\{\alpha_n:n\in\w\}$ is precisely the sequence converging to $\alpha$ used in the Fubini product. For more on these ideals and their presentations, the reader may consult \cite{Dobrinencofinalmapsonultrafilters}.

It follows directly from the definition that we can characterize the ideals $\textsf{FIN}^\cB$ via Hechler trees.

\begin{lem}\label{findisjointhechler}
    Let $\cB$ be a barrier and $X\in\fin^\cB$, then there exists a Hechler tree $H\subseteq\w^{<\w}$ such that $X\cap H\cap\cB=\emptyset$.
\end{lem}

\begin{proof}
    We prove it by induction on $\rho(\cB)$. If $\cB=[\w]^1$ and $X\in\textsf{FIN}^\cB=\textsf{FIN}$, define $H\subseteq\w^{<\w}$ such that $\emptyset\in H$, $succ_H(\emptyset)=\omega\setminus k$ where $X\subseteq k$ and $succ_H(s)=\w$ for any $s\in H\setminus\{\emptyset\}$. This $H$ clearly works.

    Let now $\cB$ be an arbitrary barrier with $\rho(\cB)>1$. We may assume that $\{k\}\notin\cB$ for every $k\in\w$. Fix $X\in\textsf{FIN}^\cB$. Thus $\{n\in\w:X\cap\cB(n)\notin\textsf{FIN}^{\cB{[n]}}\}$ is finite and we can find $k\in\omega$ that contains this set. Hence, for every $n\in\w\setminus k$, there is a Hechler tree $H_n$ with root $\{n\}$ such that $X\cap H_n\cap\cB(n)=\emptyset$. Define $H=\{\emptyset\}\cup\bigcup_{n\geq k}H_n$. It is clear that $H$ is a Hechler tree and $X\cap H\cap\cB=\emptyset$.
\end{proof}

\begin{lem}\label{lessthanbFIN}
    For every $\cX\subseteq\fin^\cB$ with $|\cX|<\b$, there is a Hechler tree $H\subseteq\w^{<\w}$ such that for every $X\in\cX$, there exists $n\in\omega$ such that 
    $$X\cap H\cap\cB\subseteq\bigcup_{i\leq n}\cB(i).$$
\end{lem}

\begin{proof}
    Let $\cX\subseteq\textsf{FIN}^\cB$ with $|\cX|<\b$. By Lemma \ref{findisjointhechler}, for every $X\in\cX$ there is a Hechler tree $H_X$ such that $X\cap H_X\cap\cB=\emptyset$. Let $f_X=f_{H_X}$ for every $X\in\cX$. Since $|\cX|<\b$ we can find $f\in\w^\w$ such that $f>^*f_X$ for every $X\in\cX$. It follows from our enumeration of $\w^{<\w}$ that if $f>^*f_X$, then there exists $n_X\in\omega$ such that if $f_X(s)>f(s)$ for some $s\neq\emptyset$ then $s\in\cB(i)$ for some $i\leq n_X$.
    Thus letting $H=H_f$ we have that $H\setminus H_X\subseteq\bigcup_{i\leq n_X}\cB(i)$ and this implies that
    \begin{equation*}
    X\cap H\cap\cB\subseteq\bigcup_{i\leq n_X}\cB(i).\qedhere\end{equation*}
    \end{proof}

The following proposition appears in \cite{kubistopologicalramsey} for the particular case of barriers of the form $[\w]^n$.  The same argument shows that this is true for every barrier, so we have the following:

\begin{propo}\label{adfamilyisBseqcomp}
    Let $\cA$ be an almost disjoint family on a countable set $N$ and let $\cB$ be a barrier. If for every function $f:\cB\to N$ there is an infinite set $X\in[\w]^\w$ such that 
    $$|\{A\in\cA:|f[\cB|X]\cap A|=\w\}|<\b$$
    then $\cF(\cA)$ is $\cB$-sequentially compact.\qed
\end{propo}

\vskip2mm

The following Lemma will also be helpful for proving that a space is not $\cB$-sequentially compact for some barrier $\cB$. We say that a barrier $\cC$ is non-trivial if $\cC\neq[\w]^1$.

\begin{lem}\label{noCseqcomp}
    Let $\cA$ be an almost disjoint family on a non-trivial barrier $\cC$, such that $|A\cap\cC(n)|\leq1$ for every $A\in\cA$ and every $n\in\omega$. If for every $E\in[\w]^\w$ there exists $A\in\cA$ such that $A\subseteq \cC|E$, then $\cF(\cA)$ is not $\cC$-sequentially compact.
\end{lem}

\begin{proof}
    As $\cC$ is non-trivial, we may assume without loss of generality that every element of the $\cC$ has size at least $2$, in particular, $\cC(m)$ is infinite for every $m\in\w$.

    Consider $i:\cC\to\cC$ to be the identity map and $E\in[\w]^\w$. We shall show that $i\rest(\cC|E)$ does not converges in $\cF(\cA)$. Since $i[\cC|E]$ is infinite, it does not converge to any isolated point in $\cC$.\\
    On the other hand, if $A\in\cA$ and $n\in\omega$, we have that $|\cC(m)\cap i[\cC|E\setminus n]|=\w$ for every $m\in E\setminus n$, but $|A\cap \cC(m)|\leq1$ implies that $i[\cC|(E\setminus n)]\nsubseteq A$ and hence it does not converge to $A$. It remains to show that it does not converge to $\infty$.\\
    Fix the neighborhood  $U=\cF(\cA)\setminus(\{A\}\cup A)$ of $\infty$ for some $A\subseteq\cC|E$ and let $n\in\omega$. As $|A\cap\cC(i)|\leq1$ for every $i\in\omega$, we can find $c\in A\subseteq\cC|E$ such that $\min(c)>n$. Thus $c\in\cC|(E\setminus n)\setminus U$ and in consequence $i[\cC|(E\setminus n)]\nsubseteq U$. Since this is true for any $n\in\omega$, we get that $i\rest(\cC|E)$ does not converge to $\infty$.
\end{proof}

The previous lemmas allow us to show that there are spaces that are $\cB$-sequentially compact but fail to be $\cC$-sequentially compact whenever their ideals associated satisfy a suitable relation, for this reason it is useful to introduce the terminology of the \Katetov order.

\begin{dfn}
    \cite{katetovproduct} Let $X$ and $Y$ be countable sets, $\cI$ and $\cJ$ ideals on $X$ and $Y$ respectively and $f:Y\rightarrow X$. 
    \begin{enumerate}
        \item $f$ is a \emph{\Katetov}function from $(Y,\cJ)$ to $(X,\cI)$ if $f^{-1}(A)\in\cJ$ for all $A\in\cI$.
        \item  $\cI\leq_K\cJ$ (\emph{$\cI$ is \Katetov below $\cJ$}) if there exists a \Katetov function from $(Y,\cJ)$ to $(X,\cI)$. 
    \end{enumerate}
\end{dfn}

%We will simple say that $f$ is a \Katetov functions when the ideals are clear by the context.
An easy consequence of the definition is the following:
\begin{lem}\label{equiv de katetov}
    Let $\cI$ and $\cJ$ ideals on countable sets $X$ and $Y$ respectively. Then the following are equivalent:
    \begin{enumerate}
        \item\label{1 de lema katetov} $\cI\not\leq_K\cJ$.
        \item\label{2 de lema katetov} For every $f:Y\rightarrow X$ there is $B\in\cJ^+$ such that $f[B]\in\cI$.
    \end{enumerate}
\end{lem}
\begin{proof}
    To see that (\ref{1 de lema katetov}) implies (\ref{2 de lema katetov}) let $f:Y\rightarrow X$. As $\cI\not\leq_K\cJ$, $f$ is not a \Katetov function, so there is $A\in\cI$ such that $B:=f^{-1}(A)\in\cJ^+$. Now $f[B]\in\cI$ since $f[B]\subseteq A\in\cI$. Conversely, to see that $\cI\not\leq_K\cJ$ let $f:Y\rightarrow X$, we want to prove that $f$ is not a \Katetov function. So let $B\in\cJ^+$ such that $A:=f[B]\in\cI$. Now $f^{-1}(A)\in\cJ^+$ since $B\subseteq f^{-1}(A)$. 
\end{proof}

Besides the ideal $\fin^\cB$, another ideal that will be useful in the constructions of counterexamples is the following:
\begin{dfn}
If $\cB$ is a barrier on $\w$, then $\mathcal{G}_c(\cB)$ is the ideal on $\cB$ such that $(\mathcal{G}_c(\cB))^+=\langle\{\cB|X: X\in[\w]^\w\}\rangle$, i.e., $S\subseteq\cB$ is an element of $\mathcal{G}_c(\cB)$ if and only if there is no $X\in[\w]^\w$ such that $\cB|X\subseteq S$.
\end{dfn}

Note that $(\cG_c([\w]^2))^+$ is the collection of all subsets of $[\w]^2$ that contain an infinite complete subgraph, that is the reason for the notation $\cG_c$. %in this sense $(\cG_c(\cB))^+$ can be interpreted as the collections of all subsets of $\cB$ that contain a complete $\cB$-hypergraph. 
Also note that the fact that $\mathcal{G}_c(\cB)$ is closed under finite unions follows from Nash-Williams Theorem.

\begin{thm}{\textnormal{($\b=\c$)}}\label{counterexampleb=c}
    Let $\cB$ and $\cC$ be two barriers such that $\fin^\cC\not\leq_K\mathcal{G}_c(\cB)$. % and assume that for every function $f:\cB\to\cC$, there exists $X\in[\w]^\w$ such that $f[\cB|X]\in\textsf{FIN}^\cC$. 
    Then there is a space that is $\cB$-sequentially compact but not $\cC$-sequentially compact.
\end{thm}

\begin{proof}
    Enumerate $\cC^\cB=\{f_\alpha:\alpha<\c\}$ and $[\w]^\w=\{E_\alpha:\alpha<\c\}$. For $E\in[\w]^\w$ let $E^\uparrow=T(\cC|E)$. Notice that for any such $E$, the tree $E^\uparrow$ is everywhere $\w$-splitting and hence it has infinite intersection with any Hechler tree $H$ inside $\cC$, that is, $|E^\uparrow\cap H\cap\cC|=\omega$. Moreover, for any $n\in E$, we have that
    $$|E^\uparrow\cap H\cap \cC(n)|=\w.$$

    Recursively construct $\{A_\alpha:\alpha<\c\}\subseteq[\cC]^\w$ and $\{X_\alpha:\alpha<\c\}\subseteq[\w]^\w$ such that for all $\beta<\alpha<\c$:
    \begin{enumerate}
        \item $f_\alpha[\cB|X_\alpha]\in\textsf{FIN}^\cC$,
        \item\label{item2counterexample} $|A_\alpha\cap A_\beta|<\w$,
        \item\label{item3counterexample} $|A_\alpha\cap f_\beta[\cB|X_\beta]|<\w$,
        \item $|A_\alpha\cap\cC(n)|\leq1$ for all $n\in\w$,
        \item\label{item5counterexample} $A_\alpha\subseteq E_\alpha^\uparrow$.
    \end{enumerate}
At step $\alpha<\c$, let $X_\alpha\in[\w]^\w$ such that % $\cB|X_\alpha\in(\cG_c(\cB))^+$ and 
$f_\alpha[\cB|X_\alpha]\in\fin^\cC$ (this set exists by Lemma \ref{equiv de katetov}). %given by the hypothesis of the Theorem (\textbf{we should clarify here that we are talking about the set that witnesses that $f_\alpha$ is not katetov}). 
Note that also, every $A_\beta$ with $\beta<\alpha$ is an element of $\textsf{FIN}^\cC$, then we can find a Hechler tree $H$ as in Lemma \ref{lessthanbFIN} such that 
$$X\cap H\cap\cC\subseteq \bigcup_{i\leq n_X}\cC(i)$$
for either $X=f_\beta[\cB|X_\beta]$ or $X=A_\beta$ for some $\beta<\alpha$. Let $A_\alpha\in[\cC]^\w$ be such that $A_\alpha\subseteq H\cap\cC\cap E_\alpha^\uparrow$ and $|A_\alpha\cap \cC(i)|\leq1$ for every $i\in\w$. It follows that $A_\alpha\cap X\subseteq \bigcup_{i\leq n_X}\cC(i)$ for some $n_X\in\w$ and as $|A_\alpha\cap\cC(i)|\leq1$ for every $i\leq n_X$, this intersection is finite and (\ref{item2counterexample}) and (\ref{item3counterexample}) hold. The remaining three properties follow from the construction.

Since $\{A\in\cA:|f_\alpha[\cB|X_\alpha]\cap A|=\w\}\subseteq\{A_\beta:\beta\leq\alpha\}$, and this set has size less that $\b$, by Proposition \ref{adfamilyisBseqcomp} we get that $\cF(\cA)$ is $\cB$-sequentially compact.

On the other hand, item (\ref{item5counterexample}) implies that $\cF(\cA)$ is not $\cC$-sequentially compact by Lemma \ref{noCseqcomp} as $\mathcal C$ has rank greater than 1.
\end{proof}

The natural question now, is when do the barriers $\cB$ and $\cC$ satisfy that $\fin^\cC\not\leq_K\mathcal{G}_c(\cB)$. %for every function $f\in\cC^\cB$ we can find an infinite set $X$ such that $f[\cB|X]\in\textsf{FIN}^\cC$. 
We will see below that this only depends on the rank of the barriers, assuming that at least the one with larger rank is uniform. See Theorem \ref{counterexamplealphaseqcomp} and Corollary
\ref{en barreras uniformes es lo mismo rango que katetov} below.

The relevance of the Kat\v{e}tov order in the study of subclasses of sequentially compact spaces is also demonstrated in the work of R. Filip\'ow, K. Kowitz and A. Kwela
\cite{filipow2023unified}, where the authors study many subclasses of sequentially compact spaces, prove some inclusions among them and find counterexamples by using the Kat\v{e}tov order relationship between some ideals naturally associated to these classes.

Recall that for $\cB$ being a barrier, $\cB[s]$ is a barrier on $\omega\setminus(\max(s)+1)$. For the ease of notation, we will do some abuse of notation: for a function $g:\cB(s)\to\cC(t)$ we can naturally associate a function $g':\cB[s]\to\cC[t]$, we use $g$ for both functions. Also, for a Hechler tree $H$, we will denote by $H[n]$ its copy above $\{n\}$, that is, a tree $T$ with root $\{n\}$ and such that $succ_T(s)=succ_H(s\setminus\{n\})$. In particular $\w^{<\w}[n]=\{s\in\w^{<\w}:n=\min(s)\}$.

\begin{thm}\label{counterexamplealphaseqcomp}
    Let $\cB,\cC$ be barriers with $\cC$ uniform and $\rho(\cB)<\rho(\cC)$. Then $\fin^\cC\not\leq_K\mathcal{G}_c(\cB)$. 
    
    %If $f:\cB\to\cC$, then we can find $X\in[\w]^\w$ and $H\subseteq\w^{<\omega}$ a Hechler tree such that $f(\cB|X)\cap H=\emptyset$.
\end{thm}

\begin{proof}
By Lemma \ref{equiv de katetov}, it is enough to prove that for every $f:\cB\to\cC$ there is $X\in[\w]^\w$ such that $f[\cB|X]\in\fin^\cC$.
The proof is by induction on $\rho(\cC)$. If $\rho(\cC)=1$ there is nothing to do, so we can assume that $\rho(\cC)>1$ and the result is true for $\cC'$ whenever $\rho(\cC')<\rho(\cC)$. We now prove it by induction on $\rho(\cB)$. %So let $f:\cB\to\cC$, we will see that $f$ is not a \Katetov function from $(\cB,\cG_c(\cB))$ to $(\cC,\fin^\cC)$. %$\fin^\cC\leq_K\mathcal{G}_c(\cB)$. 
If $\rho(\cB)=1$, we can find an infinite set $X\in[\w]^\w$ such that either $|f[\cB|X]\cap\cC(i)|\leq1$ for every $i\in\w$ or $f[\cB|X]\subseteq\cC(n)$ for some $n\in\omega$. Hence, assume that for any barriers $\cB'$ and $\cC'$ with $\rho(\cB')<\rho(\cB)$, $\rho(\cC')\leq\rho(\cC)$ and $\rho(\cB')<\rho(\cC')$ the result holds. 
We can also think in $\{\emptyset\}$ as a barrier of rank 0. In this case the result is trivial but it is worth to mention since $\cB[s]$ is a barrier of rank 0 whenever $s\in\cB$ and we will also consider these kind of barriers as part of our inductive hypothesis.
We can also assume, by possibly passing to a cofinite set, that $\rho(\cC(n))\geq \rho(\cB)$ for every $n\in\w$.

Define a function $\pi_n:\cB\to 2$ for every $n\in\omega$ given by $\pi_n(s)=0$ if and only if $f(s)\in\bigcup_{i\leq n}\cC_i$. By Nash-Williams, we can find an infinite set $Z_n\in[\w]^\w$ such that $\pi\rest(\cB|Z_n)$ is monochromatic. If for some $n\in\omega$ we can find a monochromatic set $Z_n$ with color $0$, then the Hechler tree $H$ defined such that $\emptyset\in H$, $succ_H(\emptyset)=\omega\setminus(n+1)$ and $succ_H(s)=\w$ for all $s\in H\setminus\{\emptyset\}$ works for the infinite set $X=Z_n$. Hence we can assume that there is no 0-monochromatic set for every $n\in\omega$.

Let us recursively define $\{X_n:n\in\w\}$, $\{m_i:i\in\omega\}$, $\{H_n:n\in\omega\}$ and $\{g(s,n):s\subseteq\{m_i:i<n\}\land n\in\w\}$ such that:
\begin{enumerate}
    \item $m_i=\min(X_i)$,
    \item $X_n$ is monochromatic for $\pi_n$,
    \item $X_{n+1}\subseteq X_n$,
    \item $g(s,n):\cB(s)\to \cC(n)$ is given by $g(s,n)(t)=f(t)$ if and only if $f(s)\in\cC(n)$ and otherwise let $g(s,n)(t)\in \cC(n)$ be arbitrary.
    \item $H_n\subseteq \w^{<\w}[n]$ is a Hechler tree.
    \item\label{smallimagesitem6} $g(s,n)[\cB(s)|X_n]\cap H_n=\emptyset$ for every $s\subseteq\{m_i:i\leq n\}$ with $s\neq\emptyset$.
\end{enumerate}

Since $m_i$ and $g(s,n)$ are defined from $X_n$, we only need to define $X_n$ and $H_n$. Note that if in point (\ref{smallimagesitem6}) we have that $s\notin T(\cB)$, then $\cB(s)=\emptyset$ and then the choices of $X_n$ and $H_n$ are trivial, so we will assume in the rest of the proof that $s\in T(\cB)$.

Let $X_0$ be $\pi_n$ monochromatic and let $H_0=\w^{<\w}[0]$, i.e., $H_0$ is a copy of $\w^{<\w}$ on top of $\{0\}$. 
Assume we have defined $X_i$ and $H_i$ for $i\leq n$. We can find a $\pi_{n+1}$ monochromatic set (with color 1) $X'\subseteq X_n$. 
Enumerate as $\{s_i:i\leq k\}$ for some $k\in\w$, the set $\{s\subseteq\{m_i:i\leq n\}:s\neq\emptyset\}$. 
We can recursively define 
$X_{n+1}^0\subseteq X'$, $X_{n+1}^{i+1}\subseteq X_{n+1}^i$ for $i<k$ and $H_{n+1}^i$ for $i\leq k$ such that:
\begin{itemize}
    \item[($\ast$)]\label{asterisclemma} $g(s,n)[\cB[s]|X_{n+1}^i]\cap H_{n+1}^s=\emptyset.$ 
\end{itemize}
To see that the choices of $X_{n+1}^i$ and $H_{n+1}^i$ are possible, note that $g(s,n):\cB[s]\to\cC[n]$ and $\rho(\cB[s])<\rho(\cB)\leq\rho(\cC[n])$, and then we can apply the inductive hypothesis.
Without loss of generality we can assume that $\min(X_{n+1}^i)>m_n$ for every $i\leq k$, so that no confusion arises in $(\ast)$ with the identification of $\cB(s)$ and $\cB[s]$. Also, since each $H_{n+1}^i$ is a Hechler tree for every $i\leq k$, we can define $H_{n+1}=\bigcap_{i\leq k}H_{n+1}^i[n+1]$ which is a Hechler tree on top of $\{n\}$. Define $X_{n+1}=X_{n+1}^k$. This finishes the construction.

Let $X=\{x_i:i\in\w\}$ and $H=\{\emptyset\}\cup\bigcup_{n\in\omega}H_n$. We shall show that $f[\cB|X]\cap H=\emptyset$. Fix $s\in\cB|X$, then $s=(x_{i_0},\ldots,x_{i_k})$ for some $k\in\w$ and some sequence $i_0<\cdots<i_k$. As $\{x_{i_0},\ldots,x_{i_k}\}\subseteq X_{i_0}$ and $X_{i_0}$ is $\pi_{i_0}$ monochromatic with color 1, we have that $f(s)\in C(n)$ for some $n>i_0$.
Let $s'=\{x_i\in s:i< n\}$ and notice that $s\setminus s'\subseteq X_n$ and $s'\neq\emptyset$. Hence $g(s',n)$ was considered at step $n$ in (\ref{smallimagesitem6}) and $s'$ was indexed as $s_i$ for some $i$. 
Since we have that $H\cap\w^{<\w}[n]\subseteq H_n^i$, $s\setminus s'\in\cB[s']|X_n$ and $X_n\subseteq X_n^i$, we can conclude, by $(\ast)$, that $f(s)=g(s',n)(s)\notin H$. Therefore $f[\cB|X]\cap H=\emptyset$, so $f[\cB|X]\in\fin^\cC$. 
\end{proof}

It is worth noting that the previous result gives us that for two barriers $\cB$ and $\cC$, the classes of $\cB$-sequentially compact spaces and $\cC$-sequentially compact spaces only depends on the ranks of $\cB$ and $\cC$ whenever the largest is uniform. If $\rho(\cB)<\rho(\cC)$ and $\cC$ is uniform, every $\cC$-sequentially compact space is $\cB$-sequentially compact by Corollaries \ref{alphasequentiallyequivalence} and \ref{downwardsimplicationseuqntialcompactness}. On the other hand, there are (at least consistently) examples of $\cB$-sequentially compact spaces that are not $\cC$-sequentially compact by Theorem \ref{counterexamplealphaseqcomp}. If otherwise $\rho(\cB)=\rho(\cC)$ and both of them are uniform, then the classes of $\cB$-sequentially compact and $\cC$-sequentially compact spaces coincide by Corollary \ref{alphasequentiallyequivalence}.

A more subtle analysis can be done by considering $\rho_u(\cB)$ instead of $\rho(\cB)$, but what we have done suffices to show that (consistently) the classes of $\alpha$-sequentially compact and $\beta$-sequentially compact spaces do not coincide if $\alpha\neq\beta$.

\begin{coro}{\textnormal{($\b=\c$)}}
    If $\{\alpha_n:n\in\w\}\subseteq\alpha$, there exists a space that is $\alpha_n$-sequentially compact for all $n\in\omega$ but fails to be $\alpha$-sequentially compact. In particular, for every $\alpha<\beta<\omega_1$, there exists an $\alpha$-sequentially compact space that is not $\beta$-sequentially compact. 
\end{coro}

\begin{proof}
    Fix a uniform barrier $\cB_n$ of rank $\alpha_n$ for every $n\in\omega$ and a barrier $\cC$ of rank $\alpha$. Repeat the proof of theorem \ref{counterexampleb=c} by enumerating 
    $$\bigcup_{n\in\omega}\cC^{\cB_n}=\{f_\alpha:\alpha<\c\}.$$
    The hypothesis of Theorem \ref{counterexampleb=c} for every barrier $\cB_n$ hold by Theorem \ref{counterexamplealphaseqcomp}. Finally by Corollary \ref{alphasequentiallyequivalence} we get the result.
\end{proof}

If $\cB$ and $\cC$ are barriers with the same rank and $\cB$ is uniform, then $\fin^\cC\leq_K\cG_c(\cB)$. To see this note that otherwise, by Theorem \ref{counterexampleb=c}, there is an space $X$ that is $\cB$-seq compact but not $\cC$-seq compact under $\b=\c$, but this is a contradiction to Corollary \ref{alphasequentiallyequivalence}, since the \Katetov order is absolute among Borel ideals due to Shoenfield’s absoluteness Theorem (see \cite{michaelkatetovonborel}).
We then get the following corollary:

\begin{coro}\label{en barreras uniformes es lo mismo rango que katetov}
    Let $\cB$ and $\cC$ be uniform barriers such that $\rho(\cB)\leq\rho(\cC)$. Then the following are equivalent:
    \begin{enumerate}
        \item\label{1 of equiv of katetov and rank for uniform barriers} $\fin^\cC\not\leq_K\cG_c(\cB)$.
        \item\label{2 of equiv of katetov and rank for uniform barriers} $\rho(\cB)<\rho(\cC)$.
    \end{enumerate}
\end{coro}

It is worth to point out that for (\ref{1 of equiv of katetov and rank for uniform barriers}) implies (\ref{2 of equiv of katetov and rank for uniform barriers}) we use that $\cB$ is uniform while in (\ref{2 of equiv of katetov and rank for uniform barriers}) implies (\ref{1 of equiv of katetov and rank for uniform barriers}) we use the uniformity of $\cC$.

\section{Cardinal invariants associated to barriers}\label{Cardinalinvariants}

In this section, we define and analyze several cardinal invariants associated to barriers that play an important role in the structure of $\alpha$-sequentially compact spaces. In particular, we extend a result in \cite{kubistopologicalramsey} by showing that for every $\alpha>1$, the Cantor cube $2^\kappa$ is $\alpha$-sequentially compact if and only if $\kappa<\min\{\b,\s\}$. When $\alpha=1$ it is known that $2^\kappa$ is sequentially compact if and only if $\kappa<\s$ (see \cite{vanDouwentheintegers}).

It is also shown in \cite{corralnRamsey}, that the cardinal $\parr$ is closely related to the constructions of examples of spaces that are $n$-sequentially compact but fail to be $(n+1)$-sequentially compact under $\s=\b$.
This suggest that if we aim to do similar constructions for arbitrary barriers, we should start analyzing and computing the analogous cardinal invariants for $\parr$ in arbitrary barriers.  
If $\cB$ is a barrier, $\pi:\cB\rightarrow 2$ and $H\in[\w]^\w$, then we say that $H$ is \textit{monochromatic for $\pi$}, if $\pi$ is constant on $\cB|H$. We say that $H$ is \textit{almost monochromatic} for $\pi$, if there is $F\in [H]^{<\w}$ such that $H\setminus F$ is monochromatic for $\pi$. If $\Pi$ is a family of colorings on $\cB$, then we say that $H$ is \textit{almost monochromatic for $\Pi$} if it is almost monochromatic for every $\pi\in\Pi$. We refer to $[\w]^1$ as the trivial barrier.

\begin{dfn}
    If $\cB$ is a barrier, then $\parr_\cB$ is the minimum size of a collection of colorings on $\cB$ without an almost monochromatic infinite set.
\end{dfn}

If $t\in [\w]^{<\w}\setminus\{\emptyset\}$, recall that we write $t$ as $\{t_0,\dots,t_{n-1}\}$ in the increasing way and if $I\subseteq n$, then $t\restriction I=\{t_i\mid i\in I\}$. Also note that $t\restriction 2\sqsubseteq t$ for $n\geq 2$. Given a barrier $\cB$ such that each of its elements has size at least $2$, if we let $\cB\restriction 2:=\{t\restriction 2: t\in\cB\}$ then $\cB\restriction 2=[\w]^2$.

\begin{thm}
    If $\cB$ is a non trivial barrier, then $\parr_\cB\leq\parr_2$.
\end{thm}

%\begin{proof}
 %   Let $\cS$ be a splitting family of size $\s$ and for every $S\in\cS$ let us define the coloring $\pi_S:\cB\rightarrow 2$ as $\pi_S(t)=0$ if and only if $t\subseteq S$ or $t\cap S=\emptyset$. Call $\cC=\{\pi_s\mid S\in\cS\}$.

  %  For every $S\in\cS$ there is no $1$-monochromatic infinite for $\pi_S$. Indeed, let $H\in[\w]^{\w}$ and for every $i\in 2$ let $H_i=H\cap S^i$. As $H$ is infinite there is $i\in 2$ such that $H_i$ is infinite. Now there is $s\in\cB$ such that $s\sqsubseteq H_i$, so $\pi_S(s)=0$, but $s\in \cB\restriction H$, so this shows that $H$ is not $1$-monochromatic.

    %So if $H$ is monochromatic for some $\pi_S$ then necessarily is $0$-monochromatic and $H$ can not intersect both $S^0$ and $S^1$, because if it is the case then take $a_i\in H\cap S_i$ for every $i\in 2$ and let $\widehat{H}\in [H]^\w$ such that $\{a_0,a_1\}\sqsubseteq \widehat{H}$. Now there is $t\in\cB$ such that $t\sqsubseteq\widehat{H}$, but as every element of $\cB$ has size at least $2$ then $\{a_0,a_1\}\sqsubseteq t$, so necesarly $t\cap S^i\neq\emptyset$ for both $i=0$ and $i=1$, so $\pi_S(t)=1$ and $t\in\cB\restriction H$ and this is a contradiction since $H$ was $0$-monochromatic. This implies that every $H\in [\w]^\w$ that is almost monochromatic for $\pi_S$ is such that $H\cap S=^*\emptyset$ or $H\subseteq^* S$. So if $H$ is an almost monochromatic for $\cC$ this implies that no $S\in\cS$ splits $H$, but this is not possible since $\cS$ is splitting. 
%\end{proof}

\begin{proof}
    As $\cB$ is a non trivial barrier, there is $n\in\w$ such that $\cB|(\w\setminus n)$ has all of its elements of size at least two. Without loss of generality we may assume that $n=0$. Let $\Pi=\{\pi_\alpha\mid\alpha\in\parr_2\}$ be a family of colorings on $[\w]^2$ with no almost monochromatic infinite sets. 
    Define $\widehat{\pi}_\alpha:\cB\rightarrow 2$ for every $\alpha\in\parr_2$,  given by $\widehat{\pi}_\alpha(t)=\pi_\alpha(t\restriction 2)$. 
    Note that if $H$ is a an almost monochromatic infinite set for $\widehat{\Pi}=\{\widehat{\pi}_\alpha\mid\alpha\in\parr_2\}$, then $H$ is also an almost monochromatic infinite set for $\Pi$. Hence $\widehat\Pi$ has no infinite almost monochromatic sets.
\end{proof}

Recall that if $\cB$ is a barrier,
the rank of $\cB[n]$ is strictly smaller than the rank of $\cB$ for every $n\in\w$. This allows us to show that the same idea that proves that $\parr_n=\min\{\b,\s\}$ for every $n\in\w$ in \cite{blasscardinal}, also works for arbitrary barriers.

\begin{thm}
If $\cB$ is a barrier such that each of its elements has size at least $2$, then $\parr_\cB\geq\parr_2$.    
\end{thm}
\begin{proof}
    By induction on the rank of $\cB$. Suppose that $\cB$ has rank $\gamma$ and the result is true for every barrier of rank strictly smaller than $\gamma$.
    Let $\kappa<\parr_2$ and let $\{\pi_\alpha:\alpha\in\kappa\}$ be a set of colorings on $\cB$. We are going to show that there exists a almost monochromatic infinite set for this family. 

    For every $\alpha\in\kappa$ and every $n\in\w$ let $\pi_\alpha^n=\cB[n]\rightarrow 2$ given by $\pi_\alpha^n(s)=\pi_\alpha(\{n\}\cup s)$.

    We recursively construct a collection $\{H_n: n\in\w\}\subseteq[\w]^\w$ such that:
    \begin{enumerate}
        \item $H_{n+1}\in [H_n]^\w$.
        \item $H_{n+1}$ is almost monochromatic for $\{\pi_\alpha^n\restriction(\cB[n]| H_{n}): \alpha\in\kappa\}$ 
    \end{enumerate}

    We start by defining $H_0=\w\setminus 1$ and as $B[0]$ has rank smaller than $\gamma$, by our inductive hypothesis we know that there is $H_1\in [\w\setminus 1]^\w$ such that $H_1$ is almost monochromatic for $\{\pi_\alpha^0\mid\alpha\in\kappa\}$. 
    Now suppose that we have already defined $H_n$ and consider the family $\{\pi_\alpha^n\restriction(\cB[n]|H_n): \alpha\in\kappa\}$. As $\cB[n]|H_n$ has also rank smaller than $\gamma$, there is $H_{n+1}\in[H_n]^\w$ such that $H_{n+1}$ is almost monochromatic for this family.
    
    Take $H$ an infinite pseudointersection of $\{H_n: n\in\w\}$. 
Note for every $n\in H$ and 
every $\alpha\in\kappa$ we have that $H_n$ is almost monochromatic for $\pi_\alpha^n$ with color $j(\alpha,n)\in\{0,1\}$. For every $\alpha\in\kappa$, define an increasing $f_\alpha\in H^\w$ such that $H_n\setminus f_\alpha(n)$ is monochromatic for $\pi_\alpha^n$. Also let $g\in\w^\w$ increasing such that $H\setminus g(n)\subseteq H_n$. As $\cF=\{f_\alpha:\alpha\in\kappa\}\cup \{g\}\subseteq \w^\w$ has size less than $\parr_2\leq\b$, there exists $f\in\w^\w$ strictly increasing that dominates $\cF$.

For every $i\in 2$, let $X_\alpha^i:=\{n\in H: j(\alpha,n)=i\}$. Note that $\{X_\alpha^0\mid\alpha\in\kappa\}$ is a collection of less than $\parr_2\leq\s$ subsets of $H$, so there is $X\in [H]^\w$ such that, for every $\alpha\in\kappa$ there is $j(\alpha)\in\{0,1\}$ such that $X\subseteq^* X_\alpha^{j(\alpha)}$.

 Let us define $J:=\{x_n\mid n\in\w\}\subseteq X\subseteq H$ such that $x_{n+1}>f(x_n)$  for every $n\in\w$.

     \underline{\textit{Claim:}} $J$ is almost monochromatic for $\{\pi_\alpha\mid\alpha\in\kappa\}$.

     \textit{Proof of the claim:} Let $\alpha\in\kappa$. We know that $X\subseteq^* X_\alpha^{j(\alpha)}$, $f_\alpha<^*f$ and $g<^* f$, so there exists $k\in\w$ such that:
     \begin{enumerate}
         \item $X\setminus k\subseteq X_\alpha^{j(\alpha)}$,
         \item $f_\alpha(n)<f(n)$ for every $n>k$ and
         \item $g(n)<f(n)$ for every $n>k$.
     \end{enumerate}
     Fix $s\in \cB|(J\setminus(k+1))$ and let $n=\min(s)$ and $t=s\setminus\{n\}$. 
     Then we have that $n\in X_\alpha^{j(\alpha)}$ and $$t\subseteq J\setminus f(n)\subseteq X\setminus f(n)\subseteq H\setminus f(n)\subseteq H\setminus g(n)\subseteq H_n.$$
     Notice that $\min(t)>f(n)>f_\alpha(n)$, hence $t\subseteq H_n\setminus f_\alpha(n)$. But then $\pi_\alpha(s)=\pi_\alpha^n(t)=j(\alpha,n)$ and $j(\alpha,n)=j(\alpha)$ since $n\in X_\alpha^{j(\alpha)}$. This proofs that $\pi_\alpha$ is constant on $\cB|(J\setminus(k+1))$ with color $j(\alpha)$.\hfill\ensuremath{_{Claim}\square}
     
 The proof of the previous Claim finishes the proof of the Theorem.
\end{proof}

\begin{coro}
   $\parr_\cB=\parr_2$ for every non trivial barrier $\cB$ on $\w$ and $\parr_{[\w]^1}=\s$.
\end{coro}

\begin{coro}
    If $\kappa<\parr_2$ and $\cB$ is a barrier, then $2^\kappa$ is $\cB$-sequentially compact. Moreover
    $$\parr_2=\min\{\kappa:2^\kappa\textnormal{ is not }\cB\textnormal{-sequentially compact}\}$$
    for every non trivial barrier $\cB$.
\end{coro}
\begin{proof}
    Let $f:\cB\rightarrow 2^\kappa$ for some $\kappa<\parr_2$ and for every $\alpha\in\kappa$ let $\pi_\alpha:\cB\rightarrow 2$ the $\alpha$-the projection of $f$. Thus there is $M\in[\w]^\w$ almost monochromatic with color $j(\alpha)$ for each $\alpha\in\kappa$. It follows that $f\restriction(\cB|M)$ converges to $j=\langle j(\alpha):\alpha\in\kappa\rangle\in 2^\kappa$. Conversely, if $\{\pi_\alpha:\alpha<\parr_2\}$ is a family of colorings on $\cB$ with no infinite monochromatic set, and we define $f:\cB\to 2^{\parr_2}$ by $f(b)(\alpha)=\pi_\alpha(b)$, then any infinite set $M\in[\w]^\w$ such that $f\rest(\cB|M)$ converges to a point $\langle j(\alpha):\alpha<\parr_2\rangle$, can be easily seen to be an almost monochromatic set for $\{\pi_\alpha:\alpha<\parr_2\}$, thus there is no such an infinite set $M$.
\end{proof}

The dual cardinal to $\p\a\r$ is the cardinal $\homm$. The characterization of $\p\a\r$ in terms of $\b$ and $\s$ has also a dual version using $\d$ and a variant of $\r$, the dual cardinals of $\b$ and $\s$ respectively (see \cite{blasscardinal}).

\begin{nota}
    Let $N$ be a countable set, $\cI$ an ideal on $N$ and $X,Y\subseteq N$. Then $X\subseteq_\cI Y$ means that $X\setminus Y\in\cI$.
\end{nota}

\begin{dfn}
    Given a barrier $\cB$ let, 
    \[\homm_\cB=\min\{|\cH|:\cH\subseteq[\w]^\w\ \land\ \forall S\subseteq\cB\ \exists M\in\cH\ \exists i\in 2\ (\cB|M\subseteq S^i)\}
    \]
    and
    \[\r_\cB=\min\{|\cR|:\cR\subseteq[\w]^\w\ \forall S\subseteq\cB\ \exists R\in\cR\ \exists i\in 2\ (\cB|R\subseteq_{\fin^\cB} S^i)\}.
    \]
\end{dfn}

A family as in the definition of $\homm_\cB$ is called an \emph{homogeneous family for $\cB$} or a \emph{$\cB$-homogeneous family}. Analogously a family as in the definition of $\r_\cB$ is called a \emph{reaping family for $\cB$} or a \emph{$\cB$-reaping family}. 

Recall that by the Nash-Williams theorem, for every $\pi:\cB\rightarrow 2$ there is $M\in[\w]^\w$ such that $\pi\restriction(\cB|M)$ is constant. In this way $\homm_\cB$ is well defined and $\r_\cB\leq\homm_\cB$ since every $\cB$-homogeneous family is a $\cB$-reaping family. In particular $\r_\cB$ is well defined. 

Note that if $\cB$ is the trivial barrier, i.e., $\cB=[\w]^1$, then $\r_\cB=\homm_\cB=\r$. On the other hand, if $\cB$ is a non trivial barrier, then $\homm_2\leq\homm_\cB$. To see this we can suppose that every element of $\cB$ has size at least $2$ and %let us fix $\cH$ a homegeneous family for $\cB$. Note that e
thus every coloring $\pi:[\w]^2\rightarrow 2$ induces a coloring of $\hat{\pi}:\cB\rightarrow 2$ given by $\hat{\pi}(s)=\pi(s\restriction 2)$. Now an infinite homogeneous set for $\hat{\pi}$ is also homogeneous for $\pi$, which in turns gives that any homogeneous family for $\cB$ is an homogeneous family for $[\w]^2$ and then $\max\{\d,\r_\sigma\}=\homm_2\leq\homm_\cB$. In particular $\d\leq\homm_\cB$. Ramsey ultrafilters will be helpful in our study of the cardinal $\homm_\cB$.

Recall that a non-principal ultrafilter $\cU$ on $\w$ is \emph{Ramsey} if for all $\pi:[\w]^2\rightarrow 2$, there is $U\in\cU$ such that $\pi\restriction[U]^2$ is constant.
A family $\cF\subseteq[\w]^\w$ is \emph{selective} if for every decreasing sequence $\{Y_n\mid n\in\w\}\subseteq\cF$, there is $f:\w\rightarrow\w$ such that $f[\w]\in\cF$, $f(0)\in Y_0$ and $f(n+1)\in Y_{f(n)}$ for all $n\in\w$. 
The following result, which according to D. Booth is mostly due to K. Kunen, relates Ramsey ultrafilters and selective families.

\begin{thm}\cite{boothultrafilters}
    If $\cU$ is a non principal ultrafilter, then $\cU$ is Ramsey if and only if it is selective.
\end{thm}

We will show now that a base for a Ramsey ultrafilter is an homogeneous family for every barrier $\cB$. Hence if we define $\u_R$ as the minimum character of a Ramsey ultrafilter\footnote{Let $\u_R=\c$ in case that there is no such ultrafilter.} we get that $\homm_\cB\leq\u_R$.

The next theorem is a standard result of Local Ramsey theory (see \cite{todorcevicramseyspaces} Chapter 7). It follows immediately from Theorem 7.42 in the
same way Nash-Williams' theorem follows from Ellentuck theorem. For the
convenience of reader, we sketch the argument.

\begin{thm}\label{the ramsey ultrafilters are even ramsey for barrriers}
    If $\cU$ is a Ramsey ultrafilter, $\cB$ is a barrier on $\w$ and $\pi:\cB\rightarrow 2$, then there is $U\in\cU$ such that $\pi\restriction(\cB|U)$ is constant.
\end{thm}

\begin{proof}
    By induction on the rank of the barrier. If $\cB$ has rank zero, i.e, $\cB=\{\emptyset\}$ then the result is trivial.
    
    Suppose that $\cB$ has rank $\gamma$ and the result is true for every barrier of rank smaller than $\gamma$. As $B[n]$ is a barrier on $\w\setminus(n+1)$ of rank less than $\gamma$, there is $V_n\in\cU$ such that the coloring $\pi_n:B[n]\rightarrow 2$ is constant on $\cB[n]|V_n$, where $\pi_n(s)=\pi(\{n\}\cup s)$.

    For every $n\in\w$ we know that $\pi_n\restriction(\cB[n]|V_n)$ is constant on color $j(n)\in2$. So take $V\in \cU$ and $i\in 2$ such that $j(n)=i$ for all $n\in V$.  

    Now call $U_n=(\bigcap_{m\leq n}V_n)\cap V \cap (\w\setminus n+1)$ and note that $\{U_n\mid n\in\w\}$ is a decreasing sequence of elements of $\cU$, so applying that $\cU$ is selective, there is $f\in\w^\w$ such that $f[\w]\in\cU$, $f(0)\in U_0$ and $f(n+1)\in U_{f(n)}$ for all $n\in\w$. Define $U:=f[\w]$.

    It remains to show that $\pi\restriction(\cB|U)$ is constant with value $i$. Indeed let $s=\{a_0,\dots,a_k\}\in\cB|U$. Then $a_0\in V$ and $$\{a_1,\dots,a_k\}\subseteq U_{{f(a_0)}}\subseteq U_{a_0}\subseteq V_{a_0},$$
        so that $\{a_1,\dots,a_k\}\in\cB[a_0]|V_{a_0},$. Then $\pi(s)=\pi_{a_0}(\{a_1,\dots,a_k\})=j(a_0)=i$.
\end{proof}

%\textcolor{blue}{It is well know that the forcing notion $\mathbb{U}=(\cP(\w)/\fin,\subseteq^*)$ adds a Ramsey ultrafilter $\cU$ to the ground model without adding reals (see for example \cite{halbeisen} pp. 358), thus Theorem \ref{the ramsey ultrafilters are even ramsey for barrriers} implies Theorem \ref{NashWilliamsareRamsey} in the case that the Nash-Williams is a front. To see this, let $\cF$ a front and $\pi:\cF\rightarrow 2$. Now, if we force with $\mathbb{U}$ then in $V[G]$ there is $\cU$ a Ramsey ultrafilter and $\pi\in V[G]$. Now, by Theorem \ref{the ramsey ultrafilters are even ramsey for barrriers}, there is $X\in \cU\cap V[G]$ such that  $\pi\restriction({\cF}|X)$ is constant, but as $\mathbb{U}$ does not add reals, $X\in V$ and we are done.}\carlos{If this makes sense and it is worth to leave it here, then we should change Theorem 5.9 to fronts.}

%Also note that by Theorem \ref{the ramsey ultrafilters are even ramsey for barrriers} we get that the a base for a Ramsey ultrafilter is a homogeneous family for every barrier $\cB$, this way if we define $\u_R$ as the minimum character of a Ramsey ultrafilter\footnote{Let $\u_R=\c$ in case that there is no such ultrafilter.} we get that $\homm_\cB\leq\u_R$.

\begin{nota}
If $R\in[\w]^\w$ then $\HH(R)$ is the collection of all non-empty trees $T\subseteq R^{<\w}$ such that $succ_T(s)$ is cofinite in $R$ for all $s\in T$.
\end{nota}

We will need the following folklore facts.

\begin{propo}\label{lasramasdelhechler}
    If $R\in[\w]^\w$ and $T\in\HH(R)$, then there is $A\in[R]^\w$ such that $[A]^\w\subseteq[T]$.
\end{propo}

\begin{propo}\label{cofinalidaddelfin}
$\textsf{cof}(\fin^\cB)=\d$ for every non trivial barrier $\cB$.
\end{propo}

An explicit proof for the previous result about $\textsf{cof}(\fin^\cB)$ where $\cB$ is a barrier of finite rank appears in \cite{corralnRamsey}. The same argument works for general barriers.

\begin{thm}
 $\homm_\cB=\max\{\d,\r_\cB\}$ for every non trivial barrier $\cB$ on $\w$.
\end{thm}
\begin{proof}
    We already noted that $\max\{\d,\r_\cB\}\leq\homm_\cB$, so it is enough to prove the other inequality. For this let $\cR=\{R_\beta\mid\beta \in\r_\cB\}$ be a reaping family for $\cB$ and for every $\beta\in\r_\cB$ let $\{X_\alpha^\beta\mid \alpha\in\d\}$ cofinal in $\fin^{\cB|R_\beta}$. Now by definition of $\fin^{\cB|R_\beta}$ and Proposition \ref{lasramasdelhechler}, for every $\beta\in\r_\cB$ and every $\alpha\in\d$ there are $T_\alpha^\beta\in\HH(R_\beta)$ and $A_\alpha^\beta\in [R_\beta]^\w$ such that:
    \begin{enumerate}
        \item $T_\alpha^\beta\cap X_\alpha^\beta=\emptyset$,
        \item $[A_\alpha^\beta]^\w\subseteq [T_\alpha^\beta]$.
    \end{enumerate} 

    Now let $\cH=
    \{A_\alpha^\beta\mid \beta\in\r_\beta\wedge\alpha\in\d\}$. We claim that $\cH$ is a homogeneous family for $\cB$. To see this let $S\subseteq\cB$. As $\cR$ is a reaping family for $\cB$ there is $\beta\in\r_\beta$ and $i\in 2$ such that $\cB|R_\beta\subseteq_{\fin^\cB}S^i$. Now call $X:=(\cB|R_\beta)\setminus S^i\in\fin^\cB$. Moreover, as $X\subseteq\cB|R_\beta$, then $X\in\fin^{\cB|R_\beta}$. Thus there is $\alpha\in\d$ such that $X\subseteq X_\alpha^\beta$ and consequently $X\cap T_\alpha^\beta=\emptyset$. Also we know that $[A_\alpha^\beta]^\w\subseteq[T_\alpha^\beta]$, so $X\cap (\cB|A_\alpha^\beta)=\emptyset$. Indeed, if $b\in X\cap (\cB|A_\alpha^\beta)$, then there is $Z\in[A_\alpha^\beta]^\w$ such that $b\sqsubseteq Z$ and, as $Z\in [T_\alpha^\beta]$, then $b\in T_\alpha^\beta\cap X$, which is a contradiction. This way we have that $\cB|A_\alpha^\beta\subseteq\cB|R_\beta\subseteq S^i\cup X$ and $(\cB|A_\alpha^\beta)\cap X=\emptyset$, so $\cB|A_\alpha^\beta\subseteq S^i$ and we are done.
\end{proof}

\begin{que}
    Is $\homm_\cB=\homm_2$ for every non trivial barrier $\cB$?
\end{que}

\subsection*{Acknowledgments}
This research was done while the authors were part of the Thematic Program on Set Theoretic Methods in Algebra, Dynamics and Geometry at the Fields Institute.

\bibliography{seqcomp.bib}{}
\bibliographystyle{plain}

\end{document}